\documentclass[a4paper]{siamart250106}

\usepackage{lipsum,array}
\usepackage{amsfonts}
\usepackage{graphicx}
\usepackage{epstopdf}
\usepackage{algorithmic}
\usepackage{tikz,pgf}
\usetikzlibrary{calc}
\usepackage{multicol,multirow}
\usepackage{hhline}
\usepackage{bm}
\usepackage{subcaption}
\DeclareMathOperator*{\esssup}{ess\,sup}
\ifpdf
  \DeclareGraphicsExtensions{.eps,.pdf,.png,.jpg}
\else
  \DeclareGraphicsExtensions{.eps}
\fi

\newsiamremark{remark}{Remark}
\newsiamremark{hypothesis}{Hypothesis}
\crefname{hypothesis}{Hypothesis}{Hypotheses}
\newsiamthm{claim}{Claim}

\headers{Convergence analysis of VEM for the Sobolev equation with convection}{A. Kumar, S. Kumar, and S. Yadav}

\title{Convergence analysis of virtual element methods for the Sobolev equation with convection}

\author{Ankit Kumar\thanks{Department of Mathematics, Birla Institute of Technology and Science, Pilani, Pilani Campus, Vidhya Vihar, Pilani, Rajasthan, 333031, India 
  (\email{p20210041@pilani.bits-pilani.ac.in}, \email{sangita.yadav@pilani.bits-pilani.ac.in}).}
\and Sarvesh Kumar\thanks{Corresponding author: Department of Mathematics, Indian Institute of Space Science and Technology, Thiruvananthapuram, Kerala, 699547, India
  (\email{sarvesh@iist.ac.in}).}
\and Sangita Yadav\footnotemark[1]}

\usepackage{amsopn}

\ifpdf
\hypersetup{
  pdftitle={Convergence analysis of virtual element methods for the Sobolev equation with convection},
  pdfauthor={A. Kumar, S. Kumar, and S. Yadav}
}
\fi

\usepackage[left=2cm,right=2cm,top=3cm,bottom=3cm]{geometry}
\newcolumntype{C}[1]{>{\centering\let\newline\\\arraybackslash\hspace{0pt}}m{#1}}
\begin{document}
    \newcommand{\R}{{\mathbb R}}
	\newcommand{\Z}{{\mathbb Z}}
	\newcommand{\N}{{\mathbb N}}
    \newcommand{\PS}{{\mathbb P}}
    \newcommand{\Rc}{{\mathcal R}}
    \newcommand{\eps}{\epsilon}
    \newcommand{\sg}{\sigma}
    \newcommand{\bs}{\boldsymbol}
    
\maketitle

\begin{abstract}
We explore the potential applications of virtual elements for solving the Sobolev equation with a convective term. A conforming virtual element method is employed for spatial discretization, while an implicit Euler scheme is used to approximate the time derivative. To establish the optimal rate of convergence, a novel intermediate projection operator is introduced. We discuss and analyze both the semi-discrete and fully discrete schemes, deriving optimal error estimates for both the energy norm and $L^2$-norm. Several numerical experiments are conducted to validate the theoretical findings and assess the computational efficiency of the proposed numerical methods.
\end{abstract}



\section{Introduction}

The Sobolev equations arise in modeling various physical phenomena, such as heat transfer in different media, seepage into fissured rocks, and other applications \cite{barenblatt1960basic, showalter1975sobolev, ting1974cooling}. These equations include a mixed derivative term of third order with respect to both time and space. Over the past few decades, the numerical study of the Sobolev equations has gained significant attention, leading to the development and application of various numerical methods. Finite element methods (FEM) have been widely recognized as powerful tools for solving these problems. These include the standard finite element method \cite{chen2019two, ewing1978time}, mixed finite element method \cite{dongyang2016unconditional, gao2009split, li2019expanded}, and the discontinuous Galerkin method \cite{gao2009local}, which have been employed for different formulations of the Sobolev equation. 

However, these methods face limitations in handling complex geometrical domains and irregular polygonal elements. To address these challenges, alternative approaches such as hybrid high-order (HHO) methods and the weak Galerkin finite element method (WG-FEM) have been proposed, as they can efficiently handle polygonal meshes \cite{kumar2023stabilizer, gao2017weak, xie2022hybrid}. Despite their advantages, these methods require explicit construction of local basis functions to compute the local matrices, which significantly increases the total degrees of freedom and, consequently, the computational cost.

The challenge of handling polygonal meshes without explicitly constructing local shape functions can be effectively addressed using the virtual element method (VEM), introduced by Beirão da Veiga et al. \cite{beirao2013basic} as an extension of the mimetic finite difference method. Unlike HHO and WG-FEM methods, the VEM allows for the use of non-polynomial shape functions, providing greater flexibility; see \cite{ahmad2013equivalent, beirao2016virtual}. VEM can also be considered a generalization of finite element methods on polygonal meshes. \textcolor{black}{VEM is robust on meshes with highly irregular or non-convex elements and can naturally accommodate hanging nodes without special treatment. These features make VEM a powerful tool for a wide range of applications, including parabolic problems \cite{vacca2015virtual}, elasticity \cite{da2013virtual}, fluid flow \cite{antonietti2014stream, antonietti2023virtual}, poroelasticity \cite{kumar2024numerical}, and parabolic integro-differential equations \cite{yadav2024conforming}. Moreover, VEM has been applied to partial differential equations which includes a convection term, for instance, convection-dominated-diffusion case has been studied in \cite{benedetto2016order, da2021supg}. In \cite{arrutselvi2021virtual}, the authors studied the nonlinear time-dependent convection-diffusion-reaction equation using streamlined upwind Petrov-Galerkin stabilized VEM and verified the effectiveness of the proposed method by presenting numerical results for the fully discrete scheme obtained by employing the backward Euler method.} For a detailed discussion on recent advancements in VEM, we refer readers to \cite{anaya2020virtual, brenner2017some, cangiani2017conforming}.

Recently, the VEM has been applied to the Sobolev equation as well. Zhang et al. \cite{zhang2023virtual} studied cases with constant dispersion and diffusion coefficients. Furthermore, Xu et al. \cite{xu2022conforming} and Pradhan and Deka \cite{pradhan2024optimal} extended the method to variable coefficients. The analysis in \cite{xu2022conforming} focuses on the special case where both coefficients are equal, whereas authors in \cite{pradhan2024optimal} addresses the general case, establishing optimal rate of convergence in both the energy norm and the $L^2$-norm. \textcolor{black}{These studies have primarily focused on Sobolev equations without incorporating the convection term. However, in many practical applications, such as the transport of moisture in soil, the convection term plays a fundamental role \cite{mendina2012sensitivity, wang2024influence}. In such scenarios, surface heating during the daytime increases evaporation, raises local humidity, and initiates atmospheric convection, necessitating the inclusion of convection effects in the mathematical model. The application of VEM to the Sobolev equations with convection still remains unexplored.}

The present study aims to develop the VEM for the Sobolev equation with a convection term involving variable coefficients. To construct a well-posed semi-discrete virtual element scheme for the model problem, the bilinear forms are defined following the approach in \cite{cangiani2017conforming}. Subsequently, Euler's backward difference method is employed for time discretization, leading to the corresponding fully discrete problem. The primary challenge in the analysis lies in achieving optimal convergence rate in $L^2-$norm. While the usual elliptic projection operator and the operator defined in \cite{pradhan2024optimal} ensure optimal convergence in the energy norm for our problem, but fails to achieve optimal convergence in $L^2-$norm. To address this, we define a suitable projection operator that enables optimal convergence rates in both the energy norm and $L^2$-norm. A-priori estimates in the semi-discrete situation have been established without employing Gronwall's lemma; but, in the fully discrete case, utilisation of Gronwall's lemma couldn't be avoided. The theoretical results are further validated through numerical experiments on various types of polygonal meshes.

The contents of the paper are organized as follows: In section $2$, we define our model problem and its continuous weak formulation. Section $3$ is dedicated to the construction of virtual element spaces, the semi-discrete scheme, and its well-posedness. In section $4$, we derive the error estimates for the semi-discrete problem in the energy norm and $L^2-$norm. Section $5$ presents the analysis of a fully discrete scheme. Section $6$ contains numerical tests conducted to validate the theoretical results. And section $7$ contains the concluding remarks.

In this paper, we will follow the standard notation for the functional spaces on a bounded polygonal domain $\mathcal{D}$. In particular, for any non-negative integer $s$, $|\cdot|_{s,\mathcal{D}}$ and $\|\cdot\|_{s,\mathcal{D}}$ represents the semi-norm and norm on the Sobolev spaces $H^s(\mathcal{D})$, respectively. In $L^2(\mathcal{D})$, the inner product is represented as $(\cdot, \cdot)_{\mathcal{D}}$ and the norm as $\|\cdot\|_{\mathcal{D}}$. When $\mathcal{D} = \Omega$, we will exclude the subscript $\mathcal{D}$ in the notation of inner products and norms. $\PS_s(\mathcal{D})$ will denote the space of all polynomials of degree up to $s$ on $\mathcal{D}$. Moreover, the space $L^p(0,T;H^s(\mathcal{D})), \ 1\leq p < \infty$ defined as:
\begin{equation*}
     L^p(0,T;H^s(\mathcal{D})) = \left\{v:[0,T]\to H^s(\mathcal{D}): \int_0^T\|v(\cdot,t)\|_{s,\mathcal{D}}^pdt < \infty\right\},
\end{equation*}
equipped with the norm $\|v\|_{L^p(0,T;H^s(\mathcal{D}))} = \left( \int_0^T\|v(\cdot,t)\|_{s,\mathcal{D}}^pdt\right)^{\frac{1}{p}}$. For $p = \infty$, we define the space
\begin{equation*}
    L^{\infty}(0,T;H^s(\mathcal{D})) = \left\{v:[0,T]\to H^s(\mathcal{D}):\esssup\limits_{0\leq t\leq T}\|v(\cdot,t)\|_{s,\mathcal{D}} < \infty\right\}
\end{equation*}
equipped with the norm $\|v\|_{L^{\infty}(0,T;H^s(\mathcal{D}))} = \displaystyle \esssup_{0\leq t\leq T}\|v(\cdot,t)\|_{s,\mathcal{D}}$. For any Banach space $\mathcal{W}$, the Bochner spaces $L^p(0,T;\mathcal{W})$ and $H^s(0,T;\mathcal{W})$ will be denoted by $L^p(\mathcal{W})$ and $H^s(\mathcal{W})$, respectively. Throughout this paper, we will be using $C$ as a generic positive constant which is independent of the mesh size.

\section{Model Problem}
Let $\Omega \subset \R^2$ be a polygonal domain with a Lipschitz continuous boundary $\partial \Omega$. Consider the following Sobolev equation: 
\begin{subequations}
    \begin{align}
      u_t(\bs x,t) - \nabla \cdot \left( \bs \mu(\bs x)\nabla u_t(\bs x,t) + \bs \eps(\bs x)\nabla u(\bs x,t)\right) + \bs \beta(\bs x) \cdot \nabla u(\bs x,t) + \gamma(\bs x) u(\bs x,t) &= f(\bs x,t) \quad \text{in} \ \Omega \times \mathcal{I},\label{1.mp1}\\
    u(\bs x,t) &= 0 \quad \quad \ \ \text{on} \ \partial \Omega \times \mathcal{I},\label{1.mp2}\\
    u(\bs x,0) &= u_0(\bs x) \quad \quad \quad \text{in} \ \Omega, \label{1.mp3}
    \end{align}
\end{subequations}
where $\mathcal{I} = (0,T]$, $\bs \mu(\bs x)$ and $\bs \eps(\bs x)$ are dispersion and diffusion coefficients, respectively.
\vskip5pt
\noindent\textbf{Assumptions}
\begin{itemize}
    \item[A1.] 
$\bs \mu(\bs x), \bs \eps(\bs x) \in [L^{\infty}(\Omega)]^{2\times 2}$ and are assumed to be symmetric and positive definite, i.e., there exist a constant $\underline{M} > 0,$ such that for any $\bs \xi \in \R^2, \bs \xi \neq 0$ and for almost every $\bs x \in \Omega$,
\begin{equation*}
        \bs \xi^T \bs \mu(\bs x) \bs \xi \geq \underline{M}|\bs \xi|^2, \quad \bs \xi^T \bs \eps(\bs x) \bs \xi \geq \underline{M}|\bs \xi|^2.
\end{equation*}
\item[A2.] $\gamma \in L^{\infty}(\Omega)$ and $\bs \beta \in [L^{\infty}(\Omega)]^2$ with $\nabla \cdot \bs \beta \in L^{\infty}(\Omega)$. 
\item[A3.] There exists $\sg_0 \geq 0$ such that 
\begin{equation}
    \sg(\bs x) = \gamma(\bs x) - \frac{1}{2}\nabla \cdot \bs \beta(\bs x) \geq \sg_0 \geq 0, \ \ \forall \bs x \in \Omega.
\end{equation}
\end{itemize} 

The variational problem corresponding to \eqref{1.mp1}-\eqref{1.mp3} is: we seek $u \in L^2(H_0^1(\Omega))$ such that 
\begin{equation}\label{1.wkfrm}
    \begin{split}
        m_1(u_t,v)+m_2(u_t,v)+a(u,v)+b(u,v) &= (f,v), \ \ \forall v \in H_0^1(\Omega),\\
        u(\bs x,0) &= u_0(\bs x), 
        \ \ \forall \bm{x} \in \Omega,
    \end{split}
\end{equation}
where the bilinear forms on $H_0^1(\Omega)$ are given by the following expressions:
\begin{align}\label{1.bilnrfrms}
\begin{array}{lll}
    m_1(w,v) = (w,v), &&a(w,v) = (\bs \eps \nabla w,\nabla v) + (\sg w,v),\\
    m_2(w,v) = (\bs \mu \nabla w,\nabla v),
    &&b(w,v) = \frac{1}{2}[(\bs \beta \cdot \nabla w,v) - (w,\bs \beta \cdot \nabla v)].
    \end{array}
\end{align}

These bilinear forms satisfy the following coercivity and continuity bounds. For all $w,v \in H_0^1(\Omega)$,
\begin{align}\label{1.contbds}
\begin{array}{lll}
     m_{1}(w,v) \leq \|w\|\|v\|, &&a(w,v) \leq C\|w\|_{1}\|v\|_{1},\\
    m_{2}(w,v) \leq C|w|_{1}|v|_{1},
    &&b(w,v) \leq C \|w\|_{1}\|v\|_{1}.
    \end{array}
\end{align}
Also, we have for all $w \in H_0^1(\Omega)$,
\begin{equation}\label{1.coercivitybds}
    m_{1}(w,w)=\|w\|^2, \quad m_{2}(w,w) \geq \underline{M}|w|_{1}^2, \quad a(w,w) \geq \frac{\underline{M}}{1+C_{P}}\|w\|_{1}^2, \quad b(w,w) = 0.
\end{equation}
where $C_{P}$ is the constant in the Poincaré inequality given in \cite{larson2013finite}.

Our next aim is to establish the well-posedness of the problem defined in \eqref{1.wkfrm}. In this order, we first prove the stability of the problem and then discuss the existence and uniqueness of the solution.

\begin{theorem}[Stability] Let $u$ be a solution of the problem \eqref{1.wkfrm}. Then for all $t \in \mathcal{I}$
\begin{equation}\label{1.stability_bd}
    \|u(t)\|_1^2+\int_0^t\|u(s)\|_1^2ds\leq C\left( \|u_0\|_1^2+ \int_0^t\|f(s)\|^2ds \right).
\end{equation}  
\label{1.stability_thm}
\end{theorem}

\begin{proof}
    We take $v = u$ in \eqref{1.wkfrm} to obtain
    \begin{equation*}
        m_{1}(u_{t},u) + m_{2}(u_{t},u) + a(u,u) + b(u,u) = (f,u).
    \end{equation*}
    Using of coercivity of $a(\cdot,\cdot)$ and Cauchy-Schwarz inequality, it holds that
    \begin{equation*}
        \frac{1}{2}\frac{d}{dt}m_{1}(u,u) + \frac{1}{2}\frac{d}{dt}m_{2}(u,u) + \frac{\underline{M}}{1+C_{P}}\|u\|_{1}^2 \leq \|f\|\|u\|.
    \end{equation*}
    Integrating from $0$ to $t$, and using \eqref{1.contbds} and \eqref{1.coercivitybds} we obtain
    \begin{equation*}
        \|u(t)\|^2 + |u(t)|_{1}^2 + \int_0^t\|u(s)\|_{1}^2ds \leq C\left( \|u_0\|^2 + |u_0|_{1}^2 + \int_0^t\|f(s)\|\|u(s)\| ds \right).  
    \end{equation*}
    Applying Young's inequality, for all $t \in \mathcal{I}$ we get
    \begin{equation*}
    \begin{split}
        \|u(t)\|_1^2 + \textcolor{black}{\int_0^t \|u(s)\|_{1}^2ds} \leq C \left( \|u_0\|^2 + |u_0|_{1}^2 + \int_0^t\|f(s)\|^2ds \right).
    \end{split}    
    \end{equation*}
    This shows the continuous dependence of the solution on the given data.
\end{proof}

\begin{theorem}[Existence and Uniqueness] For every $t \in \mathcal{I}$, there exists a unique solution to the problem \eqref{1.wkfrm} in $L^2(H_0^1(\Omega))$.    
\end{theorem}
\begin{proof}
    Let us consider $\{\zeta_i\}_{i=1}^{\infty}$ be a orthogonal basis for $H_0^1(\Omega)$, and suppose $\chi_m$ be a finite dimensional subspace of $H_0^1(\Omega)$ such that, $\chi_m = \text{span}\{\zeta_1, \zeta_2, \ldots, \zeta_m\}$. 
    We can write the problem: for each $t$, find $\lambda_m \in \chi_m$ such that 
    \begin{equation}\label{1.existence_wkfrm}
        m_1(\lambda_{m,t},\zeta)+m_2(\lambda_{m,t},\zeta)+a(\lambda_m,\zeta)+b(\lambda_m,\zeta) = (f,\zeta), \ \ \forall \zeta \in \chi_m,
    \end{equation}
    with $(\lambda_m(\bs x,0),\zeta) = (u_{0},\zeta)$. And hence,
    \begin{equation}\label{1.lnrcomb}
        \lambda_m = \sum_{i=1}^{m} u_{i,m}(t)\zeta_i.
    \end{equation}
    On substituting \eqref{1.lnrcomb} into \eqref{1.existence_wkfrm}, we obtain a system of linear ordinary differential equations. Hence by the theory of ordinary differential equations, there exists a solution on the interval $\mathcal{I}$.\\
    Now, on putting $\zeta = \lambda_m$ in \eqref{1.existence_wkfrm} and integrating the resulting equation from $0$ to $T$, we get the bound
    \begin{equation}
        \|\lambda_m(t)\|^2 + |\lambda_m(t)|_{1}^2 + \int_0^T\|\lambda_m(s)\|_{1}^2ds \leq C\left( \|\lambda_m(0)\|^2 + |\lambda_m(0)|_{1}^2 + \int_0^T\|f(s)\|^2ds \right),
    \end{equation}
    which is uniform in $m$. We can say that $\lambda_m$ is uniformly bounded in $L^{\infty}(H_0^1(\Omega))$ and in $L^2(H_0^1(\Omega))$. To obtain the bound for $\lambda_{m,t}$ put $\zeta = \lambda_{m,t}$ in \eqref{1.existence_wkfrm} to get
    \textcolor{black}{\begin{align*}
        m_1(\lambda_{m,t},\lambda_{m,t})+m_2(\lambda_{m,t},\lambda_{m,t})+\frac{1}{2}\frac{d}{dt}a(\lambda_m,\lambda_{m}) = (f,\lambda_{m,t})+b(\lambda_{m,t},\lambda_m), 
    \end{align*}
    using coercivity of the bilinear forms $m_1(\cdot,\cdot)$ and $m_2(\cdot,\cdot)$ and continuity of $b(\cdot,\cdot)$ we have
    \begin{align*}
    \|\lambda_{m,t}\|^2+|\lambda_{m,t}|_1^2+\frac{1}{2}\frac{d}{dt}a(\lambda_m,\lambda_{m}) \leq C\left( \|f\|\|\lambda_{m,t}\|+ \|\lambda_{m,t}\|_1\|\lambda_m\|_1 \right).
    \end{align*}
    Incorporating Young's inequality to obtain the following inequality
    \begin{align*}
        \|\lambda_{m,t}\|^2+|\lambda_{m,t}|_1^2+\frac{d}{dt}a(\lambda_m,\lambda_{m}) \leq C\left( \|f\|^2+ \|\lambda_m\|_1^2 \right),
    \end{align*}
    integrating the above equation from $0$ to $T$ and use Gronwall's lemma to get
    \begin{equation}
         \int_0^T\|\lambda_{m,t}(s)\|^2ds + \int_0^T|\lambda_{m,t}(s)|_{1}^2ds + \|\lambda_m(t)\|_{1}^2 \leq C\left( \|\lambda_m(0)\|_1^2 + \int_0^T\|f(s)\|^2ds \right).
    \end{equation}}
    Thus, $\lambda_{m,t}$ is uniformly bounded in $L^2(H_0^1(\Omega))$. By Alaoglu compactness theorem, there exists a subsequence $\{\lambda_{m_j}\}$ of $\{\lambda_m\}$ such that $\lambda_{m_j} \xrightarrow{w^*} \nu$ in $L^2(H_0^1(\Omega))$. So we have, 
    \[(\lambda_{m_j},\zeta) \xrightarrow{w^*} (\nu,\zeta),\quad (\lambda_{{m_j},t},\zeta) \xrightarrow{w^*} (\nu_t,\zeta), \quad (\nabla \lambda_{m_j},\nabla \zeta) \xrightarrow{w^*} (\nabla \nu,\nabla \zeta)\] \[\text{and }  (\nabla \lambda_{{m_j},t},\nabla \zeta) \xrightarrow{w^*} (\nabla \nu_t, \nabla \zeta).\] The limiting case of all the subsequences $\{\lambda_{m_j}\}_{j=1}^{\infty}$ and use of these weak$^{*}$ convergences in \eqref{1.existence_wkfrm} gives us
    \begin{equation}\label{1.existence_wkfrm_lim_case}
        m_1(\nu_t,\zeta)+m_2(\nu_t,\zeta)+a(\nu,\zeta)+b(\nu,\zeta) = (f,\zeta).
    \end{equation}
    Hence, by denseness of the basis $\{\zeta_i\}_{i=1}^{\infty}$ in $H_0^1(\Omega)$, we get the existence of the solution for the problem \eqref{1.wkfrm}.\\
    Now, assume that $u_1$ and $u_2$ are two distinct solutions of \eqref{1.wkfrm}, then by using \eqref{1.stability_bd}, we get that $u_1 = u_2$. Hence, the solution will be unique.
\end{proof}

\section{Virtual Element Framework}
\subsection{The virtual element space}
We shall define a polygonal or polyhedral discretization $\mathcal{T}_h$ of the domain $\Omega$ into elements $K$. For every $K \in \mathcal{T}_h$, $h_K$ denotes the diameter of $K$ where $h$ stands for mesh size of $\mathcal{T}_h$ defined as $h = \max_{K\in \mathcal{T}_h} h_K$. We use $\partial K$ to represent the boundary of $K$ and $e \subset \partial K$ for its edge. Additionally, we assume that every element $K$ satisfy the following conditions:
\begin{itemize}
    \item there exists a positive $\delta$ such that $K$ is star-shaped with respect to every point in the disk $D_{\delta}$ with radius $\delta h_K$,
    \item  $|e| \geq \delta h_K$ for any edge $e \subset \partial K$.
\end{itemize}
For every element $K \in \mathcal{T}_h$ and for $k \in \N$, the auxiliary local space of
 order $k \geq 1$ for every element $K$ is defined as follows \textcolor{black}{\cite{beirao2016virtual}}:
\begin{equation*}
    W_h^k(K) = \{ v_h \in H^1(K) \cap C^0(\partial K) : w_h|_e \in \PS_k(e) \forall e \in \partial K, \Delta w_h \in \PS_k(K)\}.
\end{equation*}
To uniquely identify a function $w_h \in  W_h^k(K)$, we can choose the following set of degrees of freedom:

\begin{enumerate}
    \item[$(\mathbb{D}_1)$] Values of $w_h$ at vertices of $K$, i.e., $w_h(V_i)$ where the $i$-th vertex of $K$ is denoted by $V_i, \ i = 1,2,\ldots ,n(K) $, 
    \item[$(\mathbb{D}_2)$] For $k>1$, for each edge $e \subset \partial K$, the values of $w_h$ at $k-1$ distinct points of $e$,
    \item[$(\mathbb{D}_3)$] For $k>1$, the internal moments of $w_h$ up to degree $k-2$, \textit{i.e.}, 
    \begin{equation*}
        \frac{1}{|K|}\int_Kw_hmd{\bs x}, \ \forall m \in \mathcal{M}_{k-2},
    \end{equation*}
\end{enumerate}
where $\mathcal{M}_k$ be the space of scaled monomials given by 
\begin{equation*}
    \mathcal{M}_k = \bigg\{\left(\frac{\bs x - \bs x_K}{h_K}\right)^{\bs \alpha}, \ |\bs \alpha| \leq k \bigg\},
\end{equation*}
for a multi-index $\bs \alpha$ such that $\bs \alpha = (\alpha_1,\alpha_2), \ |\bs \alpha| = \alpha_1 + \alpha_2$ and $\bs x^{ \bs \alpha} = x_1^{\alpha_1}x_2^{\alpha_2}$ and $\bs x_K$ is the barycenter of $K$. Note that $\mathcal{M}_k$ will form a basis for the space $\PS_k(K)$. In Figure \ref{fig:1.dofs}, the degrees of freedom for $k = 1,2$ and $3$ for a polygonal element are demonstrated. The red circles represent degrees of freedom on vertices and edges, and the blue squares represent internal degrees of freedom.

\begin{center}
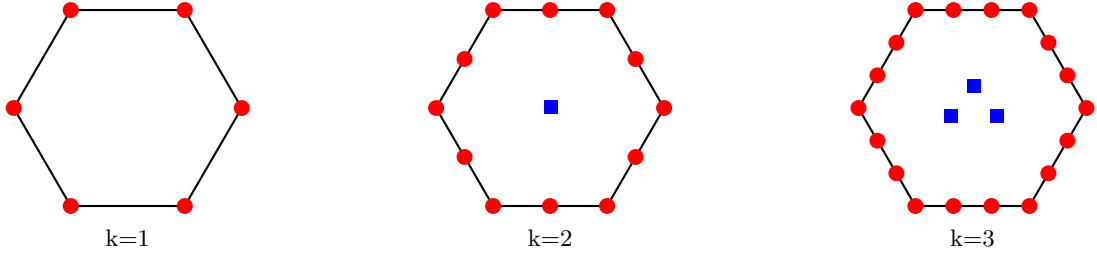
\begin{figure}[H]
    \begin{minipage}{0.32\textwidth}
    \centering
        \begin{tikzpicture}
            \foreach \i in {0,60,120,180,240,300} {
                \coordinate (V\i) at (\i:1.5);
            }
            \draw[thick] (V0) -- (V60) -- (V120) -- (V180) -- (V240) -- (V300) -- cycle;
            
            \foreach \i in {0,60,120,180,240,300} {
                \fill[red] (V\i) circle (3pt);
            }
        \end{tikzpicture}\\
        {\small k=1}
   \end{minipage}
    \begin{minipage}{0.32\textwidth}
        \centering
        \begin{tikzpicture}
    \def\r{1.5} 
    \foreach \i in {0,60,120,180,240,300} {
        \coordinate (V\i) at (\i:\r);
    }
    \draw[thick] (V0) -- (V60) -- (V120) -- (V180) -- (V240) -- (V300) -- cycle;
    \foreach \i in {0,60,120,180,240,300} {
        \fill[red] (V\i) circle (3pt);
    }
    \foreach \i/\j in {0/60, 60/120, 120/180, 180/240, 240/300, 300/0} {
        \coordinate (M\i\j) at ($(V\i)!0.5!(V\j)$);
        \fill[red] (M\i\j) circle (3pt);
    }
    \coordinate (C) at (-0.07,-0.07);
    \fill[blue] (C) rectangle ++(5pt,5pt);
\end{tikzpicture}\\
       {\small k=2}
   \end{minipage}
    \begin{minipage}{0.32\textwidth}
    \centering
        \begin{tikzpicture}
    \def\r{1.5} 
    \foreach \i in {0,60,120,180,240,300} {
        \coordinate (V\i) at (\i:\r);
    }
    \draw[thick] (V0) -- (V60) -- (V120) -- (V180) -- (V240) -- (V300) -- cycle;
    \foreach \i in {0,60,120,180,240,300} {
        \fill[red] (V\i) circle (3pt);
    }
    \foreach \i/\j in {0/60, 60/120, 120/180, 180/240, 240/300, 300/0} {
        \coordinate (M1\i\j) at ($(V\i)!1/3!(V\j)$);
        \coordinate (M2\i\j) at ($(V\i)!2/3!(V\j)$);
        \fill[red] (M1\i\j) circle (3pt);
        \fill[red] (M2\i\j) circle (3pt);
    }
    \coordinate (C1) at (-0.07,0.2);
    \coordinate (C2) at (-0.37,-0.2);
    \coordinate (C3) at (0.23,-0.2);
    \fill[blue] (C1) rectangle ++(5pt,5pt);
    \fill[blue] (C2) rectangle ++(5pt,5pt);
    \fill[blue] (C3) rectangle ++(5pt,5pt);
\end{tikzpicture}\\
        {\small k=3}
   \end{minipage}
    \caption{Degrees of freedom for a polygonal element}
    \label{fig:1.dofs}
\end{figure}
\end{center}

Next, we introduce $\Pi_{k}^{\nabla}: W_h^k(K) \to \PS_k(K)$ as the elliptic projection operator defined by
\begin{align*}
    (\nabla \Pi_{k}^{\nabla}w_h,\nabla p_k)_K &= (\nabla w_h, \nabla p_k)_K, \ \forall p_k \in \PS_k(K),\\
    P_0(\Pi_{k}^{\nabla}w_h) &= P_0w_h,
\end{align*}
where 
\begin{align*}
    P_0w_h = \left\{
    \begin{array}{cc}
        \displaystyle\frac{1}{n(K)}\sum\limits_{i=1}^{n(K)}w_h(V_i), & k=1 \\
      \frac{1}{|K|}\int_{K}w_h~d\bs x, & k>1.
    \end{array} \right.
\end{align*}

In order to compute $L^2-$projection, the local virtual space is defined as 
\begin{equation*}
    W_h(K) = \left\{ w_h \in W_h^k(K) : \int_K \Pi_{k}^{\nabla}w_hp~d\bs x = \int_K w_hp~d\bs x \ \forall p \in \PS_k(K)/\PS_{k-2}(K) \right\},
\end{equation*}
where $\PS_k(K)/\PS_{k-2}(K)$ is the space of polynomials of degree less than or equal to $k$ on $K$ that are $L^2-$orthogonal to $\PS_{k-2}(K)$. Additionally, we shall consider the $L^2-$projection operator $\Pi_k^0: L^2(K) \to \PS_k(K)$ defined by:
\begin{equation*}
    (\Pi_k^0w_h,p_k)_K = (w_h,p_k)_K, \ \forall p_k \in \PS_k(K).
\end{equation*}
Next, the global equivalent of the discrete space of order $k$ can be constructed as
\begin{equation*}
    W_h = \{ w_h \in H_0^1(\Omega) : w_h|_K \in W_h(K), \forall K \in \mathcal{T}_h \}.
\end{equation*}

\subsection{Discrete bilinear forms}
We now define the discrete bilinear forms to approximate $m_1^K(\cdot,\cdot)$, $m_2^K(\cdot,\cdot)$, $a^K(\cdot,\cdot)$, $b^K(\cdot,\cdot)$ which are restrictions of bilinear forms defined in \eqref{1.bilnrfrms} to the element $K$. To make the discrete local bilinear forms computable, we make use of 
the $L^2-$projection operator $\Pi_k^0$. For every $w_h,v_h \in W_h(K)$ the local discrete bilinear forms are defined by:
\begin{equation}\label{1.dscrtbilnrfrms}
    \begin{split}
        m_{1,h}^K(w_h,v_h) &= (\Pi_k^0w_h,\Pi_k^0v_h)_K + S_1^K((I-\Pi_k^0)w_h,(I-\Pi_k^0)v_h), \\
    m_{2,h}^K(w_h,v_h) &= (\bs \mu \Pi_{k-1}^0 \nabla w_h,\Pi_{k-1}^0\nabla v_h)_K+S_2^K((I-\Pi_k^0)w_h,(I-\Pi_k^0)v_h),\\
    a_h^K(w_h,v_h) &= (\bs\eps \Pi_{k-1}^0 \nabla w_h,\Pi_{k-1}^0 \nabla v_h)_K + (\sg\Pi_k^0w_h,\Pi_k^0v_h)_K + S_3^K((I-\Pi_k^0)w_h,(I-\Pi_k^0)v_h),\\
    b_h^K(w_h,v_h) &= \frac{1}{2}[(\bs \beta \cdot \Pi_{k-1}^0 \nabla w_h,\Pi_{k}^0v_h)_K - (\Pi_{k}^0w_h,\bs \beta \cdot \Pi_{k-1}^0\nabla v_h)_K],
    \end{split}
\end{equation}
where the stabilisation terms $S_1^K(\cdot,\cdot), S_2^K(\cdot,\cdot)$ and $S_3^K(\cdot,\cdot)$ are defined as:
\begin{equation*}
    \begin{split}
        S_1^K(w_h,v_h) &= |K|\sum_{i=1}^{N^{K}} \text{dof}_i(w_h)\text{dof}_i(v_h), \ \ w_h,v_h \in ker(\Pi_k^0), \\
        S_2^K(w_h,v_h) &= {\mu}_K\sum_{i=1}^{N^{K}} \text{dof}_i(w_h)\text{dof}_i(v_h), \ \ w_h,v_h \in ker(\Pi_k^0),\\
        S_3^K(w_h,v_h) &= \left({\eps}_K+\sg_K|K| \right)\sum_{i=1}^{N^{K}} \text{dof}_i(w_h)\text{dof}_i(v_h), \ \ w_h,v_h \in ker(\Pi_k^0),
    \end{split}
\end{equation*}
where $\mu_K,\eps_K$ and $\sg_K > 0$ are constant approximations of $\bs \mu, \bs \eps$ and $\sg$ on $K$ respectively, and $N^K$ denotes the number of degrees of freedom of $W_h(K)$. More details on the construction of  $S_1^K(\cdot,\cdot), S_2^K(\cdot,\cdot)$ and $S_3^K(\cdot,\cdot)$ can be found in \cite{cangiani2017conforming, zhang2023virtual}. Note that these stabilization terms are symmetric and positive definite such that for all \textcolor{black}{$w_h \in W_h(K) \cap ker(\Pi_k^0)$} there exist positive constants $\underline{\alpha}_i,\overline{\alpha}_i, i = 1,2,3$ independent of $K$ and $h_K$ satisfying 
\begin{equation}\label{1.stabbds}
    \begin{split}
        \underline{\alpha}_1\|w_h\|_{K}^2 &\leq S_1^K(w_h,w_h) \leq \overline{\alpha}_1\|w_h\|_{K}^2,\\
        \underline{\alpha}_2m_2^{K}(w_h,w_h) &\leq S_2^K(w_h,w_h) \leq \overline{\alpha}_2m_2^{K}(w_h,w_h),\\
        \underline{\alpha}_3a^{K}(w_h,w_h) &\leq S_3^K(w_h,w_h) \leq \overline{\alpha}_3a^{K}(w_h,w_h).
    \end{split}
\end{equation}
For all $p_k \in \PS_k(K)$ and for all $w_h \in W_h(K)$,
\begin{align}\label{1.polycons}
m_{1,h}^K(p_k,w_h) = m_1^K(p_k,w_h).
\end{align}
The global counterparts of the discrete bilinear forms can be assembled as:
\begin{align*}
    \begin{array}{lll}
    m_{1,h}(w_h,v_h) = \displaystyle \sum_{K \in \mathcal{T}_h} m_{1,h}^K(w_h,v_h), &&a_{h}(w_h,v_h) = \displaystyle \sum_{K \in \mathcal{T}_h} a_{h}^K(w_h,v_h),\\
    m_{2,h}(w_h,v_h) = \displaystyle \sum_{K \in \mathcal{T}_h} m_{2,h}^K(w_h,v_h),
    &&b_{h}(w_h,v_h) = \displaystyle \sum_{K \in \mathcal{T}_h} b_{h}^K(w_h,v_h).
    \end{array}
\end{align*}
We now state the continuity and coercivity bounds for the discrete bilinear forms. For all $w_h,v_h \in W_h$, we have the following bounds:
\begin{equation}\label{1.contybdsdscrt}
   \hspace{-1em} \begin{array}{ll}
    m_{1,h}(w_h,v_h) \leq (1+\overline{\alpha}_1)\|w_h\|\|v_h\|, &
    m_{2,h}(w_h,v_h) \leq (1+\overline{\alpha}_2)\overline{M}|w_h|_{1
    }|v_h|_{1},\\
    a_h(w_h,v_h) \leq \max\{(1+\overline{\alpha}_3)\overline{M},\|\sg\|_{\infty}\}\|w_h\|_{1}\|v_h\|_{1},
    &b_h(w_h,v_h) \leq \|\bs \beta\|_{\infty} \|w_h\|_{1}\|v_h\|_{1},
    \end{array}
\end{equation}
and for all $w_h \in W_h$,
\begin{equation}\label{1.coercivitybdsdscrt}
    m_{1,h}(w_h,w_h) \geq \min\{1,\underline{\alpha}_1\}\|w_h\|^2, \quad
    m_{2,h}(w_h,w_h) \geq \min\{1,\underline{\alpha}_2\}\underline{M}|w_h|_{1}^2,
    \quad b_h(w_h,w_h) = 0.
\end{equation}
It also holds that for all $w_h \in W_h$ \cite{cangiani2017conforming},
\begin{equation}\label{1.coercivitybdfa}
   a_h(w_h,w_h) \geq C\|w_h\|_{1}^2.
\end{equation}
Next, we define the projection operator $P^{\nabla}: H_0^1(\Omega) \to W_h$ as 
 \begin{equation*}
     m_{2,h}(P^{\nabla}w,v_h) = m_2(w,v_h) \quad \forall v_h \in W_h,
 \end{equation*}
 and its estimates are given by the following lemma (for proof see \cite{pradhan2024optimal}):

 \begin{lemma}
     Let $u \in H_0^1(\Omega)\cap H^{k+1}(\Omega)$. Then there exists a constant $C$ such that the following bounds hold:
     \begin{align}
      \|P^{\nabla}u-u\|+h   |P^{\nabla}u-u|_{1} &\leq Ch^{k+1}|u|_{k+1} \label{1.elliptic_proj_est_h1}.
     \end{align}
 \end{lemma}

\subsection{The semi-discrete problem}
We are now ready to state the semi-discrete version of our problem in the virtual element space defined by $W_h$. The virtual element approximate problem of \eqref{1.wkfrm} is: find $u_h(t) \in W_h$ for each $t\in\mathcal{I}$ such that for all $v_h\in W_h$
\begin{equation}\label{1.semidscrtprblm}
\begin{split}
    m_{1,h}(u_{h,t}(t),v_h) + m_{2,h}(u_{h,t}(t),v_h) + a_h(u_h(t),v_h) + b_h(u_h(t),v_h) &= (f_h,v_h),\\
    u_h(\cdot,0) &= P^{\nabla}u_0,
\end{split}
\end{equation}
where ${f_h}_{|_K} = \Pi_k^0f$ for all $K \in \mathcal{T}_h$.

We next establish the well-posedness of the above problem. We will first show that the semi-discrete scheme is stable, and then we will prove the existence of a unique solution. To get the continuous dependence result for $u_h$, we follow the procedure of Theorem \ref{1.stability_thm} for \eqref{1.semidscrtprblm} and obtain for all $t \in \mathcal{I}$,
\begin{equation}\label{1.stability}
     \|u_h(t)\|_1^2 + \int_0^t\|u_h(s)\|_{1}^2ds \leq C\left( \|u_h(0)\|_1^2+\int_0^t\|f_h(s)\|^2ds  \right).
\end{equation}  

If $\displaystyle u_h(t) = \sum_{i=1}^{N^W}P_i(t)\phi_i$, where $\{\phi_i\}_{i=1}^{N^W}$ be the finite-dimensional basis of the space $W_h$ and $N^W$ is the dimension of $W_h$, then \eqref{1.semidscrtprblm} is reduced to the following system of first-order ODEs:
\begin{equation}\label{1.semidiscrtode}
\begin{split}
    G \frac{dP(t)}{dt} + H P(t) &= F(t),\\
    P(0) &= P_0,
\end{split}
\end{equation}
where $P_0$ is known, $P(t) = (P_1(t), P_2(t), \ldots, P_{N^W}(t))^T$, $G$ and $H$ are $N^W\times N^W$ symmetric matrices defined as $G_{ij} = m_{1,h}(\phi_i,\phi_j)+m_{2,h}(\phi_i,\phi_j)$ and  $H_{ij} = a_h(\phi_i,\phi_j)+b_h(\phi_i,\phi_j)$, F denotes the vector $(F_1,F_2,\ldots , F_{N^W})$ defined by $F_j = (f_h,\phi_j)$. The coercivity of the bilinear forms $m_{1,h}(\cdot,\cdot)$ and $m_{2,h}(\cdot,\cdot)$ implies that $G$ is positive-definite. By the theory of ODE, there exists a unique solution to the problem \eqref{1.semidiscrtode}.

\section{Error Estimates}
In this section, we drive the optimal error estimates for the proposed semi-discrete scheme. Prior to that, we present some useful results. Let $w_I \in W_h$ be the function that interpolates $w$ in $W_h$ as $\text{dof}_i(w_I) = \text{dof}_i(w), i = 1,2,\ldots ,N^W$. Then the following estimates hold (see \cite{vacca2015virtual})
\begin{equation}\label{1.interpestimate}
    \|w - w_I\|_K + h_K|w-w_I|_{1,K} \leq Ch_K^{k+1}|w|_{k+1,K}.
\end{equation}

\textcolor{black}{We have the following results on the consistency of the bilinear forms. For proof we refer to [Theorem 6.3, \cite{cangiani2017conforming}].
\begin{lemma}\label{1.lemma_consis_err}
    Assume $u \in H^{k+1}(\Omega)$. Then for all $K \in \mathcal{T}_h$ and for all $v_h \in V_h(K)$, it holds that
    \begin{align}
        m_{2}^K(\Pi_k^0 u,v_h) - m_{2,h}^K(\Pi_k^0 u,v_h) &\leq Ch_K^k|u|_{k+1,K}\|v_h\|_{1,K}, \label{1.constncy_err1}\\
        a^K(\Pi_k^0 u,v_h) - a_h^K(\Pi_k^0 u,v_h) &\leq Ch_K^k|u|_{k+1,K}\|v_h\|_{1,K}, \label{1.constncy_err2}\\
        b^K(\Pi_k^0 u,v_h) - b_h^K(\Pi_k^0 u,v_h) &\leq Ch_K^k|u|_{k+1,K}\|v_h\|_{1,K}. \label{1.constncy_err3}
    \end{align}
\end{lemma}}

Next, we consider two bilinear forms $\mathcal{A}: H_0^1(\Omega) \times H_0^1(\Omega) \to \R$ and $\mathcal{A}_h: W_h \times W_h \to \R$ given by
\begin{align*}
    \mathcal{A}(w,v) &= m_2(w_t,v) + a(w,v) + b(w,v),\\
    \mathcal{A}_h(w_h,v_h) &= m_{2,h}(w_{h,t},v_h) + a_h(w_h,v_h) + b_h(w_h,v_h).
\end{align*}
\subsection{Intermediate Projection}
To facilitate the analysis of establishing the optimal order of convergence, we define a new intermediate projection operator $\Rc^h: H_0^1(\Omega) \to W_h$ such that for all $w_h \in W_h$
\begin{equation}\label{1.interproj}
    \mathcal{A}_h(\Rc^hw,v_h) = \mathcal{A}(w,v_h), \quad \forall t \in \mathcal{I},
\end{equation}
with $\Rc^hw(\cdot,0) = P^{\nabla}w(\cdot,0)$. The approximation properties of the intermediate projection operator are established in the subsequent theorem.
\begin{theorem}
    Let $u$ be the solution of equation \eqref{1.wkfrm} and assume $u,u_t \in {L^2(H^{k+1}(\Omega))}$. Then 
    \begin{align}
        |\Rc^hu-u|_{1} &\leq Ch^{k} \left( \|u\|_{{L^2(H^{k+1}(\Omega))}} + |u|_{k+1} \right), \label{1.interprojerrh1}\\
        \|\Rc^hu-u\| &\leq Ch^{k+1}\left(\|u\|_{L^2(H^{k+1}(\Omega))} + \|u_t\|_{L^2(H^{k+1}(\Omega))}+|u|_{k+1}\right). \label{1.interprojerrl2}
    \end{align}
\end{theorem}

\begin{proof}
    We begin the proof by setting $\rho_h = \Rc^hu - P^{\nabla}u$. Then it is evident that
    \textcolor{black}{\begin{align*}
        \mathcal{A}_h(\rho_h,\rho_h) &= \mathcal{A}_h(\Rc^hu,\rho_h)-\mathcal{A}_h(P^{\nabla}u,\rho_h) \\
        &= \mathcal{A}_h(\Rc^hu,\rho_h)-\sum\limits_{K\in \mathcal{T}_h}\left( m_{2,h}^K(P^{\nabla}u_t,\rho_h) + a_h^K(P^{\nabla}u,\rho_h) + b_h^K(P^{\nabla}u,\rho_h)\right) \\
        &= \mathcal{A}(u,\rho_h) -\sum\limits_{K\in \mathcal{T}_h}\left( m_{2}^K(u_{t},\rho_h) + a_h^K(P^{\nabla}u-\Pi_k^0 u,\rho_h) + b_h^K(P^{\nabla}u-\Pi_k^0 u,\rho_h)\right)\\
        &\quad -\sum\limits_{K\in \mathcal{T}_h}\left( a_h^K(\Pi_k^0 u,\rho_h) + b_h^K(\Pi_k^0 u,\rho_h)\right) \\
        &= \sum\limits_{K\in \mathcal{T}_h}\left( a^K(u-\Pi_k^0 u,\rho_h) + b^K(u-\Pi_k^0 u,\rho_h) - a_h^K(P^{\nabla}u-\Pi_k^0 u,\rho_h) - b_h^K(P^{\nabla}u-\Pi_k^0 u,\rho_h)\right) \\
        &\quad +\sum\limits_{K\in \mathcal{T}_h}\left( a^K(\Pi_k^0 u,\rho_h) - a_h^K(\Pi_k^0 u,\rho_h) + b^K(\Pi_k^0 u,\rho_h) - b_h^K(\Pi_k^0 u,\rho_h)\right).
    \end{align*}}
    \textcolor{black}{By coercivity of $a_h(\cdot,\cdot)$ and Lemma \ref{1.lemma_consis_err}, we have}
    \begin{align*}
        \frac{1}{2}\frac{d}{dt}m_{2,h}(\rho_h,\rho_h) + \|\rho_h\|_{1}^2 
        &\leq \textcolor{black}{C \sum\limits_{K\in \mathcal{T}_h}\left( \|u - \Pi_k^0 u\|_{1,K} + \|P^{\nabla}u - u\|_{1,K}\right) \|\rho_h\|_{1,K} }\\
        &\textcolor{black}{\quad +\sum\limits_{K\in \mathcal{T}_h}\left( a^K(\Pi_k^0 u,\rho_h) - a_h^K(\Pi_k^0 u,\rho_h) + b^K(\Pi_k^0 u,\rho_h) - b_h^K(\Pi_k^0 u,\rho_h)\right)} \\
        &\leq C \sum\limits_{K\in \mathcal{T}_h} h_K^k |u|_{k+1,K} \|\rho_h\|_{1,K} \quad \quad \left(\text{using } \eqref{1.elliptic_proj_est_h1} \right)\\
        &\leq C h^k |u|_{k+1} \|\rho_h\|_{1}.
    \end{align*}
    Using Young's inequality, we achieve
    \begin{equation*}
        \frac{1}{2}\frac{d}{dt}m_{2,h}(\rho_h,\rho_h) + \|\rho_h\|_{1}^2 \leq C\frac{h^{2k}}{2} |u|_{k+1}^2 + \frac{1}{2}\|\rho_h\|_{1}^2.
    \end{equation*}
   It yields
    \begin{equation*}
         \frac{d}{dt}m_{2,h}(\rho_h,\rho_h) + \|\rho_h\|_{1}^2 \leq Ch^{2k} |u|_{k+1}^2.
    \end{equation*}
    Integrating from $0$ to $t$ and making use of \eqref{1.contybdsdscrt} and \eqref{1.coercivitybdsdscrt} for $m_{2,h}(\cdot,\cdot)$ we obtain
    \begin{equation*}
       \min\{1,\underline{\alpha}_2\}\underline{M}|\rho_h(t)|_{1}^2 + \int_0^t\|\rho_h(s)\|_{1}^2ds \leq Ch^{2k}\int_0^t |u(s)|_{k+1}^2ds + (1+\overline{\alpha}_2)\overline{M}|\rho_h(0)|_{1}^2,
    \end{equation*}
    recall that  $\Rc^hu(\bs x,0) = P^{\nabla}u(\bs x,0)$, hence
    \begin{equation*}
       \min\{1,\underline{\alpha}_2\}\underline{M}|\rho_h(t)|_{1}^2 + \int_0^t\|\rho_h(s)\|_{1}^2ds
       \leq Ch^{2k}\int_0^t |u(s)|_{k+1}^2ds,
    \end{equation*}
 Thus
    \begin{equation}\label{1.priorprojenergy}
       |\rho_h(t)|_{1}^2+ \int_0^t\|\rho_h(s)\|_{1}^2ds \leq Ch^{2k} \|u\|_{L^2(H^{k+1}(\Omega))}^2,
    \end{equation}
    Note that
    \begin{equation*}
        |\Rc^hu-u|_{1} \leq  |\Rc^hu-P^{\nabla}u|_{1} + |P^{\nabla}u-u|_{1}.
    \end{equation*}
    Using \eqref{1.priorprojenergy} and \eqref{1.interpestimate} we obtain
    \begin{equation}\label{1.projerrorenergy}
        |\Rc^hu-u|_{1} \leq Ch^{k} \left( \|u\|_{L^2(H^{k+1}(\Omega))} + |u|_{k+1} \right),
    \end{equation}
    and
    \begin{equation}\label{1.projerrorenergy2}
        \int_0^t\|\Rc^hu-u\|_{1}^2ds \leq Ch^{2k} \|u\|_{L^2(H^{k+1}(\Omega))}^2.
    \end{equation}
    To derive the estimate in $L^2-$norm, we define a dual problem: For given $t\in\mathcal{I}$, find $\psi:[0,t) \to H^2(\Omega) \bigcap H_0^1(\Omega)$ such that 
    \begin{equation}\label{1.dualprblm}
        -\nabla \cdot \left( \bs \eps \nabla \psi - \bs \mu \nabla \psi_s\right) - \bs \beta \cdot \nabla \psi + (\gamma -\nabla \cdot \bs \beta)\psi = \rho_h, \ \ s<t
    \end{equation}
    with $\psi(t) = 0$. 
    Then $\Omega$ being convex makes $\psi$ to satisfy the following  bound
    \begin{equation}\label{1.regestimate}
        \|\psi\|_{H^1(H^2(\Omega))} \leq C\|\rho_h\|_{L^2(L^2(\Omega))}.
    \end{equation}
    On multiplying equation \eqref{1.dualprblm} by $\rho_{h,s} = \Rc^hu_s-P^{\nabla}u_s$ we get
    \begin{align*}
        (\rho_{h,s},\rho_h) &= (\rho_{h,s},-\nabla \cdot \left( \bs \eps \nabla \psi - \bs \mu \nabla \psi_s\right) - \bs \beta \cdot \nabla \psi + (\gamma - \nabla \cdot \bs \beta)\psi)\\
        &= (\nabla \rho_{h,s},\bs \eps \nabla \psi)-(\nabla \rho_{h,s},\bs \mu \nabla \psi_s)-(\rho_{h,s},\bs \beta \cdot \nabla \psi)+(\rho_{h,s},(\gamma-\nabla \cdot \bs \beta)\psi)\\
        &= (\nabla \rho_{h,s},\bs \eps \nabla \psi)-(\nabla \rho_{h,s},\bs \mu \nabla \psi_s)+(\rho_{h,s},\sg \psi)-\frac{1}{2}(\rho_{h,s},\bs \beta \cdot \nabla \psi)+\frac{1}{2}(\bs \beta \cdot \nabla \rho_{h,s},\psi).
    \end{align*}
    We can rewrite the above equation as
    \begin{align*}
        \frac{1}{2}\frac{d}{ds}(\rho_h,\rho_h) &= \frac{d}{ds}(\nabla \rho_h,\bs \eps \nabla \psi)-(\nabla \rho_h,\bs \eps \nabla \psi_s) + \frac{d}{ds}(\rho_h,\sg \psi)-(\nabla \rho_{h,s},\bs \mu \nabla \psi_s)-(\rho_h,\sg \psi_s)\\
        &~~-\frac{1}{2}\frac{d}{ds}(\rho_h,\bs \beta \cdot \nabla \psi)+\frac{1}{2}(\rho_h,\bs \beta \cdot \nabla \psi_s)+\frac{1}{2}\frac{d}{ds}(\bs \beta \cdot \nabla \rho_h,\psi)-\frac{1}{2}(\bs \beta \cdot \nabla\rho_h,\psi_s).
    \end{align*}
    Integrating from $0$ to $t$ and using $\psi(t) = 0$ and $\Rc^hu_0 = P^{\nabla}u(0)$, we arrive at
    \begin{align*}
        \frac{1}{2}\|\rho_h(t)\|^2 &= \int_0^t\left( -(\nabla \rho_h,\bs \eps \nabla \psi_s)-(\nabla \rho_{h,s},\bs \mu \nabla \psi_s)-(\rho_h,\sg \psi_s)+\frac{1}{2}(\rho_h,\bs \beta \cdot \nabla \psi_s)\right.\\
        &\left.~~-\frac{1}{2}(\bs \beta \cdot \nabla\rho_h,\psi_s)\right)ds \\
        &= \int_0^t\left( -(\bs \eps \nabla \rho_h,\nabla \psi_s)-(\bs \mu \nabla \rho_{h,s},\nabla \psi_s)-(\sg \rho_h,\psi_s)+\frac{1}{2}(\rho_h,\bs \beta \cdot \nabla \psi_s)\right.\\
        &\left.~~-\frac{1}{2}(\bs \beta \cdot \nabla \rho_h,\psi_s)\right)ds,
    \end{align*}
    using definition of $\mathcal{A}(\cdot, \cdot)$, we can write  
    \begin{align*}
        \frac{1}{2}\|\rho_h(t)\|^2 &= -\int_0^t \mathcal{A}(\rho_h,\psi_s)ds \\
        &= \int_0^t (-\mathcal{A}(\rho_h,\psi_s-\psi_{sI})-\mathcal{A}(\Rc^hu,\psi_{sI})+\mathcal{A}(P^{\nabla}u,\psi_{sI}))ds,
    \end{align*}
    it yields
    \begin{equation}\label{1.projl2erreqn}
        \begin{split}
            \frac{1}{2}\|\rho_h(t)\|^2 = \int_0^t \left(-\mathcal{A}(\rho_h,\psi_s-\psi_{sI})+\mathcal{A}_h(\Rc^hu,\psi_{sI})-\mathcal{A}(\Rc^hu,\psi_{sI}) +\mathcal{A}(P^{\nabla}u,\psi_{sI})-\mathcal{A}(u,\psi_{sI})\right)ds.
        \end{split}
    \end{equation}
    Now, using the continuity of the bilinear forms, we have 
    \begin{align*}
        \mathcal{A}(\rho_h,\psi_s-\psi_{sI}) &= \sum\limits_{K\in \mathcal{T}_h}\left(m_2^K(\rho_{h,s},\psi_s-\psi_{sI}) + a^K(\rho_h,\psi_s-\psi_{sI})+b^K(\rho_h,\psi_s-\psi_{sI}) \right) \\
        &\leq \sum\limits_{K\in \mathcal{T}_h}\left( \overline{M}|\rho_{h,s}|_{1,K}|\psi_s-\psi_{sI}|_{1,K}+\|\bs \beta\|_{\infty}\|\rho_h\|_{1,K}\|\psi_s-\psi_{sI}\|_{1,K}  \right. \\
        & \quad \left.+\max\{\overline{M},\|\sg\|_{\infty}\}\|\rho_h\|_{1,K}\|\psi_s-\psi_{sI}\|_{1,K} \right)\\
        &\leq C \sum\limits_{K\in \mathcal{T}_h}\left( |\rho_{h,s}|_{1,K} + \|\rho_h\|_{1,K} \right)\|\psi_s-\psi_{sI}\|_{1,K}\\
        &\leq C \sum\limits_{K\in \mathcal{T}_h}h_K\left( |\rho_{h,s}|_{1,K} + \|\rho_h\|_{1,K} \right)|\psi_s|_{2,K}
    \end{align*}
    \textit{i.e,}
    \begin{equation}\label{1.projapprox1}
        \mathcal{A}(\rho_h,\psi_s-\psi_{sI}) \leq Ch\left( |\rho_{h,s}|_{1} + \|\rho_h\|_{1} \right)|\psi_s|_{2}.
    \end{equation}    
    \textcolor{black}{We write
    \begin{equation}\label{1.projapprox21}
        \begin{split}
            \mathcal{A}_h(\Rc^hu,\psi_{sI})-\mathcal{A}(\Rc^hu,\psi_{sI}) &= \sum\limits_{K\in \mathcal{T}_h} \left( \mathcal{A}_h^K(\Rc^hu-\Pi_k^0u,\psi_{sI}-\Pi_1^0\psi_s)-\mathcal{A}^K(\Rc^hu-\Pi_k^0u,\psi_{sI}-\Pi_1^0\psi_s) \right.\\
            &~~+ \mathcal{A}_h^K(\Rc^hu-\Pi_k^0u,\Pi_1^0\psi_s) - \mathcal{A}^K(\Rc^hu-\Pi_k^0u,\Pi_1^0\psi_s)\\
            &~~+\left. \mathcal{A}_h^K(\Pi_k^0u,\psi_{sI}) -\mathcal{A}^K(\Pi_k^0u,\psi_{sI})\right).
        \end{split}
    \end{equation}
    Using \eqref{1.interpestimate} we derive that
    \begin{align*}
       \mathcal{A}_h^K(\Rc^hu-\Pi_k^0u,\psi_{sI}-&\Pi_1^0\psi_s)-\mathcal{A}^K(\Rc^hu-\Pi_k^0u,\psi_{sI}-\Pi_1^0\psi_s) \\
       &= m_{2,h}^K(\Rc^hu_s-\Pi_k^0 u_s,\psi_{sI}-\Pi_1^0\psi_s) + a_h^K(\Rc^hu-\Pi_k^0 u,\psi_{sI}-\Pi_1^0\psi_s) \\
        & ~~+ b_h^K(\Rc^hu-\Pi_k^0 u,\psi_{sI}-\Pi_1^0\psi_s)-m_{2}^K(\Rc^hu_s-\Pi_k^0 u_s,\psi_{sI}-\Pi_1^0\psi_s)\\
        & ~~-a^K(\Rc^hu-\Pi_k^0 u,\psi_{sI}-\Pi_1^0\psi_s)-b^K(\Rc^hu-\Pi_k^0 u,\psi_{sI}-\Pi_1^0\psi_s) \\
        &\leq C\left( h_K^{k+1}\left(|u_s|_{k+1,K}+|u|_{k+1,K}\right)+ h_K\left(|\Rc^hu_s-u_s|_{1,K}+\|\Rc^hu-u\|_{1,K}\right) \right)|{\psi_s}|_{2,K}.
    \end{align*}
    Now, 
    \begin{align*}
        \mathcal{A}_h^K(\Pi_k^0u,\psi_{sI}) -\mathcal{A}^K(\Pi_k^0u,\psi_{sI}) &= m_{2,h}^K(\Pi_k^0u_s,\psi_{sI}) -m_2^K(\Pi_k^0u_s,\psi_{sI}) + a_{h}^K(\Pi_k^0u,\psi_{sI}) - a^K(\Pi_k^0u,\psi_{sI})\\
        &~~+ b_{h}^K(\Pi_k^0u,\psi_{sI}) - b^K(\Pi_k^0u,\psi_{sI}).
    \end{align*}
    Again using Lemma \ref{1.lemma_consis_err}, we deduce that
    \begin{align*}
        \mathcal{A}_h^K(\Pi_k^0u,\psi_{sI}) -\mathcal{A}^K(\Pi_k^0u,\psi_{sI}) \leq Ch^{k+1}_K\left(|u|_{k+1,K}+|u_s|_{k+1,K}\right)|\psi_s|_{2,K},
    \end{align*}
    and 
    \begin{align*}
         \mathcal{A}_h^K(\Rc^hu-\Pi_k^0u,\Pi_1^0\psi_s) &- \mathcal{A}^K(\Rc^hu-\Pi_k^0u,\Pi_1^0\psi_s) \\
         &= m_{2,h}^K(\Rc^hu_s-\Pi_k^0u_s,\Pi_1^0\psi_s) - m_2^K(\Rc^hu_s-\Pi_k^0u_s,\Pi_1^0\psi_s) + a_{h}^K(\Rc^hu-\Pi_k^0u,\Pi_1^0\psi_s)\\
         &~~- a^K(\Rc^hu-\Pi_k^0u,\Pi_1^0\psi_s) + b_{h}^K(\Rc^hu-\Pi_k^0u,\Pi_1^0\psi_s) - b^K(\Rc^hu-\Pi_k^0u,\Pi_1^0\psi_s) \\
         &\leq Ch_K^{k+1}\left(|u_s|_{k+1,K}+|u|_{k+1,K}\right)|\psi_s|_{2,K}.
    \end{align*}
    On putting values in \eqref{1.projapprox21}, we obtain
    \begin{equation}\label{1.projapprox2}
        \mathcal{A}_h(\Rc^hu,\psi_{sI})-\mathcal{A}(\Rc^hu,\psi_{sI}) \leq C\left(h^{k+1}\left( |u_s|_{k+1} + |u|_{k+1}\right) + h\left(|\Rc^h u_s-u_s|_1+\|\Rc^hu-u\|_1 \right)\right)|\psi_s|_2.
    \end{equation}}
    Now,
    \begin{align*}
        \mathcal{A}(P^{\nabla}u,\psi_{sI})-\mathcal{A}(u,\psi_{sI}) &= \mathcal{A}(P^{\nabla}u-u,\psi_{sI}-\psi_s)+\mathcal{A}(P^{\nabla}u-u,\psi_s)\\
        &\leq \|P^{\nabla}u-u\|_{1,h}\|\psi_{sI}-\psi_s\|_{1,h}+m_{2}(P^{\nabla}u_s-u_s,\psi_s)\\
        &~~+a(P^{\nabla}u-u,\psi_s)+b(P^{\nabla}u-u,\psi_s)\\
        &\leq \|P^{\nabla}u-u\|_{1,h}\|\psi_{sI}-\psi_s\|_{1,h}+(P^{\nabla}u_s-u_s,-\nabla \cdot (\mu \nabla \psi_s))\\
        &~~+(P^{\nabla}u-u,-\nabla \cdot (\eps \nabla \psi_s))+(P^{\nabla}u-u,\sg \psi_s)\\
        &~~+\frac{1}{2}(P^{\nabla}u-u,-\nabla \cdot (\bs \beta \psi_s))+\frac{1}{2}(P^{\nabla}u-u,-\bs \beta \cdot \nabla \psi_s))\\
        &\leq \|P^{\nabla}u-u\|_{1}\|\psi_{sI}-\psi_s\|_{1}+(\|P^{\nabla}u_s-u_s\|+\|P^{\nabla}u-u\|)\|\psi_s\|_{2},
    \end{align*}
    using \eqref{1.interpestimate}, we get
    \begin{equation}\label{1.projapprox3}
        \mathcal{A}(P^{\nabla}u,\psi_{sI})-\mathcal{A}(u,\psi_{sI}) \leq Ch^{k+1}\left( |u|_{k+1} + |u_s|_{k+1} \right) \|\psi_s\|_{2}.
    \end{equation}
  Substitute \eqref{1.projapprox1},\eqref{1.projapprox2} and \eqref{1.projapprox3} in \eqref{1.projl2erreqn} and obatin
  \begin{align*}
 \|\rho_h(t)\|^2 &= C\int_0^t\left[h^{k+1}\left(|u|_{k+1} + |u_s|_{k+1}\right)+h\left(|\Rc^hu_s-u_s|_{1}+\|\Rc^hu-u\|_{1}\right)\right]\|{\psi_s}\|_{2}ds,
     \end{align*}
     Use of H\"{o}lder's inequality, \eqref{1.projerrorenergy} and \eqref{1.projerrorenergy2} yield
     \begin{align*}
         \|\rho_h(t)\|^2 
         \leq Ch^{k+1}\left(\|u\|_{L^2(H^{k+1}(\Omega))} + \|u_t\|_{L^2(H^{k+1}(\Omega))}\right)\|\psi\|_{H^1(H^{2}(\Omega))} .
     \end{align*}
      The regularity estimate given by \eqref{1.regestimate}  helps us to achieve
     \begin{equation*}
    \|\rho_h(t)\|^2  \leq Ch^{k+1}\left(\|u\|_{L^2(H^{k+1}(\Omega))} + \|u_t\|_{L^2(H^{k+1}(\Omega))}\right)\|\rho_h\|_{L^2(L^{2}(\Omega))}.   
     \end{equation*}
     Thus, we arrive at
     \begin{equation}\label{1.projl2estimatepre}
         \|\rho_h\|_{L^\infty(\mathcal{I};L^{2}(\Omega))}^2  \leq Ch^{k+1}\left(\|u\|_{L^2(H^{k+1}(\Omega))} + \|u_t\|_{L^2(H^{k+1}(\Omega))}\right)\|\rho_h\|_{L^2(L^{2}(\Omega))}.
     \end{equation}
     
    \noindent As we know that
     \begin{equation*}
         \|\rho_h\|_{L^2(L^{2}(\Omega))} \leq C\|\rho_h\|_{L^\infty(\mathcal{I};L^{2}(\Omega))}.
     \end{equation*}
     Using the above bound in equation \eqref{1.projl2estimatepre}, we get
     \begin{equation*}
         \|\Rc^hu(t)-P^{\nabla}u(t)\|_{L^\infty(\mathcal{I};L^{2}(\Omega))} \leq C h^{k+1}\left(\|u\|_{L^2(H^{k+1}(\Omega))} + \|u_t\|_{L^2(H^{k+1}(\Omega))}\right).
     \end{equation*}
     Making use of \eqref{1.elliptic_proj_est_h1} we obtain, for each $t \in \mathcal{I}$
     \begin{equation}\label{1.interprojl2estimate}
         \|\Rc^hu(t)-u(t)\| \leq C h^{k+1}\left(\|u\|_{L^2(H^{k+1}(\Omega))} + \|u_t\|_{L^2(H^{k+1}(\Omega))}+|u|_{k+1}\right).
     \end{equation}
     It completes the proof.
\end{proof}
\subsection{A priori Error Estimates}
We next establish the optimal order of convergence in both $L^2$ norm and $H^1$ norm.
\begin{theorem}
    Let $u$ and $u_h$ be the solutions of equations \eqref{1.wkfrm} and \eqref{1.semidscrtprblm}, respectively. Assume that $f,u,u_t,u_{tt} \in {L^2(H^{k+1}(\Omega))}$ and $u_0 \in H^{k+1}(\Omega)$. Then,
    \begin{equation*}\label{1.finalerrestimatel2}
    \begin{split}
        \|u_h-u\| &+ h|u_h-u|_1 \\
        &\leq Ch^{k+1}\left( \|f\|_{L^2(H^{k+1}(\Omega))} + \|u\|_{L^2(H^{k+1}(\Omega))} + \|u_t\|_{L^2(H^{k+1}(\Omega))} + \|u_{tt}\|_{L^2(H^{k+1}(\Omega))}+
        |u_0|_{k+1} \right).
        \end{split}
    \end{equation*}
\end{theorem}

\begin{proof}
    Set $u_h-u = (u_h-\Rc^hu)+(\Rc^hu-u) = \theta + \rho$. \\
    We now deduce the estimates for $\theta$ and estimates for $\rho$ we have already established in \eqref{1.interprojerrh1} and \eqref{1.interprojerrl2}. Then, from \eqref{1.wkfrm} and \eqref{1.semidscrtprblm} we have
    \begin{align*}
        m_{1,h}(\theta_{t},v_h) + m_{2,h}(\theta_{t},v_h) + a_h(\theta,v_h) + b_h(\theta,v_h) &= (f_h,v_h) - m_{1,h}(\Rc^hu_{t},v_h)- m_{2,h}(\Rc^hu_{t},v_h) \\
        &~~- a_h(\Rc^hu,v_h) - b_h(\Rc^hu,v_h) \\
        &= (f_h,v_h) - m_{1,h}(\Rc^hu_{t},v_h) - \mathcal{A}_h(\Rc^hu,v_h)\\
        &= (f_h,v_h) - m_{1,h}(\Rc^hu_{t},v_h) - \mathcal{A}(u,v_h),
    \end{align*}
    it yields,
    \begin{align}
        m_{1,h}(\theta_{t},v_h) + m_{2,h}(\theta_{t},v_h) + a_h(\theta,v_h) + b_h(\theta,v_h) &= (f_h-f,v_h)+ m_1(u_{t},v_h) - m_{1,h}(\Rc^hu_{t},v_h) \label{1.semidscrterreqn} \\
        & = \kappa_1 + \kappa_2. \nonumber
    \end{align}
    The estimates for $\kappa_1$ and $\kappa_2$ will be proved the same as \cite{vacca2015virtual}. But since we have defined a novel intermediate projection, it is better to present it. Now,
    \begin{align*}
        \kappa_1 &= (f_h-f,v_h) = \sum\limits_{K\in \mathcal{T}_h} (\Pi_k^0f-f,v_h)_K \\
        &\leq \sum\limits_{K\in \mathcal{T}_h} \|\Pi_k^0f-f\|_K\|v_h\|_K\\
        &\leq \sum\limits_{K\in \mathcal{T}_h} Ch_K^{k+1}|f|_{k+1,K}\|v_h\|_K,
    \end{align*}
    it implies,
    \begin{equation}\label{1.estimateT1}
        \kappa_1 \leq  Ch^{k+1}|f|_{k+1}\|v_h\|,
    \end{equation}
    by the polynomial consistency \eqref{1.polycons} and continuity bounds \eqref{1.contbds} and \eqref{1.contybdsdscrt}, we have
    \begin{align*}
        \kappa_2 &= m_1(u_{t},v_h) - m_{1,h}(\Rc^hu_{t},v_h)\\
        &= \sum\limits_{K\in \mathcal{T}_h} \left(m_1^K(u_{t},v_h)-m_1^K(\Pi_k^0u_{t},v_h) + m_{1,h}^K(\Pi_k^0u_{t},v_h) - m_{1,h}^K(\Rc^hu_{t},v_h)\right) \\
        &= \sum\limits_{K\in \mathcal{T}_h} \left(m_1^K(u_{t}-\Pi_k^0u_{t},v_h) + m_{1,h}^K(\Pi_k^0u_{t}-\Rc^hu_{t},v_h) \right)\\
        &\leq \sum\limits_{K\in \mathcal{T}_h} \left( \|u_{t}-\Pi_k^0u_{t}\|_K\|v_h\|_K + (1+\overline{\alpha}_1)\|\Pi_k^0u_{t}-\Rc^hu_{t}\|_K\|v_h\|_K \right)\\
        &\leq C\sum\limits_{K\in \mathcal{T}_h} \left( \|u_{t}-\Pi_k^0u_{t}\|_K + \|u_{t}-\Rc^hu_{t}\|_K\right) \|v_h\|_K,
    \end{align*}
    so,
    \begin{equation}\label{1.estimateT2}
        \kappa_2 \leq Ch^{k+1} \left( \|u_t\|_{L^2(H^{k+1}(\Omega))} + \|u_{tt}\|_{L^2(H^{k+1}(\Omega))}+|u_t|_{k+1} \right)  \|v_h\|.
    \end{equation}
    On putting values of \eqref{1.estimateT1} and \eqref{1.estimateT2} in \eqref{1.semidscrterreqn} and taking $v_h = \theta$, we achieve
    \begin{align*}
        m_{1,h}(\theta_{t},\theta) + m_{2,h}(\theta_{t},\theta) + a_h(\theta,\theta) &+ b_h(\theta,\theta)\\
        &\leq Ch^{k+1} \left(  |u_t|_{k+1} + |f|_{k+1} +  \|u_t\|_{L^2(H^{k+1}(\Omega))} + \|u_{tt}\|_{L^2(H^{k+1}(\Omega))} \right) \|\theta\|,
    \end{align*}
    using coercivity of $a_h(\cdot,\cdot)$ and Young's inequality, it yields that
     \begin{align*}
        \frac{1}{2}\frac{d}{dt}m_{1,h}(\theta,\theta)+\frac{1}{2}\frac{d}{dt}m_{2,h}(\theta,\theta) &+\|\theta\|_{1}^2\\
        &\leq Ch^{2(k+1)}\left( \|u_t\|_{L^2(H^{k+1}(\Omega))}^2 + \|u_{tt}\|_{L^2(H^{k+1}(\Omega))}^2 +|u_t|_{k+1}^2 + |f|_{k+1}^2 \right) + \frac{1}{2}\|\theta\|^2,
    \end{align*}
    Observe that the $H^1$-norm in the third term on the left side allows us to combine $\|\theta\|$ with this term. There will be no necessity to employ Gronwall's inequality in this context. Now, integrate from $0$ to $t$ and recall that $\theta(\bs x, 0) = 0$, we get
    \begin{align*}
        m_{1,h}(\theta(t),\theta(t))+m_{2,h}(\theta(t),\theta(t))+  2\int_0^t\|\theta(s)\|_{1}^2 &\leq Ch^{2(k+1)}\int_0^t\left(\|u_t(s)\|_{L^2(H^{k+1}(\Omega))}^2  + \|u_{tt}(s)\|_{L^2(H^{k+1}(\Omega))}^2\right.\\
        &~~\left. + |u_t(s)|_{k+1}^2 + |f(s)|_{k+1}^2 \right)ds, 
    \end{align*}
    it implies that 
    \begin{equation*}\label{1.errestimatesprior}
        \begin{split}
            \min\{1,\underline{\alpha}_1\}\|\theta(t)\|^2+\min\{1,\underline{\alpha}_2\}\underline{M}|\theta(t)|_{1}^2 &\leq Ch^{2(k+1)}\left(\|u_t\|_{L^2(H^{k+1}(\Omega))}^2 + \|u_{tt}\|_{L^2(H^{k+1}(\Omega))}^2 +\|f\|_{L^2(H^{k+1}(\Omega))}^2 \right)
        \end{split}
    \end{equation*}
    from the above equation, it is easy to say that,
    \begin{align}\label{1.errestimatebothnorm}
        \|\theta(t)\| + |\theta(t)|_1 \leq Ch^{k+1}\left(\|u_t\|_{L^2(H^{k+1}(\Omega))} + \|u_{tt}\|_{L^2(H^{k+1}(\Omega))} + \|f\|_{L^2(H^{k+1}(\Omega))} \right)
    \end{align}
    Note the fact that 
    \begin{align}
        |u|_{k+1} = \left|\int_0^tu_t(\cdot,s)ds + u_0\right|_{k+1} \leq \|u_t\|_{L^1(H^{k+1}(\Omega))}+|u_0|_{k+1} \leq \|u_t\|_{L^2(H^{k+1}(\Omega))} + |u_0|_{k+1}. \label{1.rhobound}
    \end{align}
    Using triangle inequality, \eqref{1.interprojerrh1}-\eqref{1.interprojerrl2} and \eqref{1.errestimatebothnorm}-\eqref{1.rhobound}, we obtain
    \begin{equation*}
    \begin{split}
        \|u_h-u\| &+ h|u_h-u|_1 \\
        &\leq Ch^{k+1}\left( \|f\|_{L^2(H^{k+1}(\Omega))} + \|u\|_{L^2(H^{k+1}(\Omega))} + \|u_t\|_{L^2(H^{k+1}(\Omega))} + \|u_{tt}\|_{L^2(H^{k+1}(\Omega))} + |u_0|_{k+1} \right).
    \end{split}
    \end{equation*}
    Here, we complete the proof.
\end{proof}

\section{The Fully Discrete Scheme}
In this section, we employ the Euler backward method for the temporal derivatives. For $n = 0,1,\ldots,N,~t_n = n\tau$ where $\tau = T/N$. We denote $f^j = f(\cdot,t_j),~u_t^j = u_t(\cdot,t_j)$. The fully discrete scheme to the problem \eqref{1.mp1}-\eqref{1.mp3} is defined as: For $n \geq 1$ we seek $\displaystyle U^n \in W_h$ such that for all $v_h \in W_h$,
\begin{equation}\label{1.fullydscrtprblm}
    \begin{split}
        m_{1,h}\left(\frac{\displaystyle U^n-U^{n-1}}{\tau},v_h\right)+m_{2,h}\left(\frac{\displaystyle U^n-U^{n-1}}{\tau},v_h\right)+a_h(U^n,v_h)+b_h(U^n,v_h)
        &= (f_h^n,v_h),\\
        U^0 &= \Rc^hu(\bs x,0),
    \end{split}
\end{equation}
where $U^n \approx u_h(\cdot,t_n)$ and $f_h^n = f_h(\cdot,t_n)$.
\begin{theorem}
    Let $u^n$ be the solution of the problem \eqref{1.wkfrm} at $t = t_n$ and $\displaystyle U^n$ be the solution of the fully discrete problem \eqref{1.fullydscrtprblm}. Assume that $f \in \displaystyle L^{\infty}(H^{k+1}(\Omega))$ with $u,u_t \in {L^2(H^{k+1}(\Omega))}$ and $u_{tt} \in \displaystyle L^1(L^{2}(\Omega))$. Then,
    \begin{equation}\label{1.fullydscrterrestimatel2}
        \begin{split}
           \|U^n-u^n\| + h|U^n-u^n|_{1} &\leq Ch^{k+1} \left( \|f\|_{L^{\infty}(H^{k+1}(\Omega))} + \|u\|_{L^2(H^{k+1}(\Omega))} + \|u_t\|_{L^2(H^{k+1}(\Omega))} + |u_0|_{k+1} \right)\\
           &~~+ C\tau \|u_{tt}\|_{L^1(L^2(\Omega))},
        \end{split}
    \end{equation}
\end{theorem}

\begin{proof}
    Set $U^j-u^j = \left( U^j-\Rc^hu^j \right) + \left( \Rc^hu^j-u^j \right) = \theta^j - \rho^j, ~j = 1, 2\ldots,N$.\\
    From \eqref{1.wkfrm} and \eqref{1.fullydscrtprblm}, we have
    \begin{align*}
        m_{1,h}\left(\frac{\theta^j-\theta^{j-1}}{\tau},v_h\right)&+m_{2,h}\left(\frac{\theta^j-\theta^{j-1}}{\tau},v_h\right)+a_h(\theta^j,v_h)+b_h(\theta^j,v_h) = (f_h^j,v_h)-a_h(\Rc^hu^j,v_h)\\
        &-b_h(\Rc^hu^j,v_h)-m_{1,h}\left(\frac{\Rc^hu^j-\Rc^hu^{j-1}}{\tau},v_h\right)-m_{2,h}\left(\frac{\Rc^hu^j-\Rc^hu^{j-1}}{\tau},v_h\right),
    \end{align*}
    using the definition of the intermediate projection, we obtain the following equation:
    \begin{equation}\label{1.erreqnfullydscrt}
        \begin{split}
             m_{1,h}\left(\frac{\theta^j-\theta^{j-1}}{\tau},v_h\right)&+m_{2,h}\left(\frac{\theta^j-\theta^{j-1}}{\tau},v_h\right)+a_h(\theta^j,v_h)+b_h(\theta^j,v_h) \\
             &= (f_h^j-f^j,v_h) + m_1(u_t^j,v_h)-m_{1,h}\left(\frac{\Rc^hu^j-\Rc^hu^{j-1}}{\tau},v_h\right)\\
             &= L_1+L_2,
        \end{split}
    \end{equation}
    In order to determine the bounds for $L_1$ and $L_2$, we first bound $L_1$ using approximation property of $\Pi_k^0$ such that
    \begin{align}
        L_1 &= (f_h^j-f^j,v_h) = \sum\limits_{K\in \mathcal{T}_h} (\Pi_k^0f^j-f^j,v_h)_K \nonumber \\
        &\leq \sum\limits_{K\in \mathcal{T}_h}\|\Pi_k^0f^j-f^j\|_K\|v_h\|_K \leq \sum\limits_{K\in \mathcal{T}_h}Ch_K^{k+1}|f^j|_{k+1,K}\|v_h\|_K \nonumber \\
        &\leq Ch^{k+1}|f^j|_{k+1}\|v_h\|. \label{1.fullydscrtestimate1}
    \end{align}
    And
    \begin{align*}
        L_2 &= m_1(u_t^j,v_h)-m_{1,h}\left(\frac{\Rc^hu^j-\Rc^hu^{j-1}}{\tau},v_h\right) \\
        &= \sum\limits_{K\in \mathcal{T}_h} \left( m_1^K(u_t^j,v_h)-m_{1,h}^K\left(\frac{\Rc^hu^j-\Rc^hu^{j-1}}{\tau},v_h\right) \right) \\
        &= \sum\limits_{K\in \mathcal{T}_h} \left( m_1^K\left(u_t^j-\frac{u^j-u^{j-1}}{\tau},v_h\right)+m_1^K\left(\frac{u^j-u^{j-1}}{\tau}-\frac{\Pi_k^0(u^j-u^{j-1})}{\tau},v_h\right) \right) \\
        &~~+\sum\limits_{K\in \mathcal{T}_h} m_{1,h}^K\left(\frac{\Pi_k^0(u^j-u^{j-1})}{\tau}-\frac{\Rc^hu^j-\Rc^hu^{j-1}}{\tau},v_h\right) \\
        &\leq C\sum\limits_{K\in \mathcal{T}_h} \left( \|u_t^j-\frac{u^j-u^{j-1}}{\tau}\|_K + \frac{1}{\tau}\|(u^j-u^{j-1})-\Pi_k^0(u^j-u^{j-1})\|_K \right) \|v_h\|_K \\
        &~~+C\sum\limits_{K\in \mathcal{T}_h}\left( \frac{1}{\tau}\|\Pi_k^0(u^j-u^{j-1})-\Rc^h(u^j-u^{j-1})\|_K \right)\|v_h\|_K \\
        &\leq \frac{C}{\tau}\sum\limits_{K\in \mathcal{T}_h} \left( \|(u^j-u^{j-1})-\Pi_k^0(u^j-u^{j-1})\|_K +\|(u^j-u^{j-1})-\Rc^h(u^j-u^{j-1})\|_K \right. \\
        &~~ \left. +\|\tau u_t^j-(u^j-u^{j-1})\|_K \right)\|v_h\|_K \\
        &\leq \frac{C}{\tau}\sum\limits_{K\in \mathcal{T}_h} \left( \|(u^j-u^{j-1})-\Rc^h(u^j-u^{j-1})\|_K+\|\tau u_t^j-(u^j-u^{j-1})\|_K \right.\\
        &~~\left.+h_K^{k+1}|u^j-u^{j-1}|_{k+1,K}\right)\|v_h\|_K,
    \end{align*}
    it implies,
    \begin{equation}\label{1.fullydscrtestimate2}
    \begin{split}
        L_2 \leq \frac{C}{\tau}\left( \|(u^j-u^{j-1})-\Rc^h(u^j-u^{j-1})\|+\|\tau u_t^j-(u^j-u^{j-1})\|+h^{k+1}|u^j-u^{j-1}|_{k+1}\right)\|v_h\|.
    \end{split}
    \end{equation}
    On putting values from \eqref{1.fullydscrtestimate1} and \eqref{1.fullydscrtestimate2} to \eqref{1.erreqnfullydscrt} and taking $v_h = \theta^j$, then \eqref{1.rhobound} yields that
    \begin{equation}\label{1.fullydscrterrbound1}
        \begin{split}
             m_{1,h}\left(\frac{\theta^j-\theta^{j-1}}{\tau},\theta^j\right)&+m_{2,h}\left(\frac{\theta^j-\theta^{j-1}}{\tau},\theta^j\right)+a_h(\theta^j,\theta^j) \\
             &\leq \frac{Ch^{k+1}}{\tau} \left( \|u\|_{L^2(H^{k+1}(\Omega))} + \|u_{t}\|_{L^2(H^{k+1}(\Omega))} + |u_0|_{k+1} + \tau|f^j|_{k+1} \right) \|\theta^j\|\\
             &~~+\frac{C}{\tau}\left(h^{k+1}|u^j-u^{j-1}|_{k+1} + \|\tau u_t^j-(u^j-u^{j-1})\| \right) \|\theta^j\|,
    \end{split}
    \end{equation}
    We now follow the same process as \cite{zhang2023virtual} and observe that
    \begin{align*}
        m_{1,h}\left(\frac{\theta^j-\theta^{j-1}}{\tau},\theta^j\right) = \frac{1}{\tau}\left( \frac{1}{2} m_{1,h}(\theta^j-\theta^{j-1},\theta^j-\theta^{j-1})+\frac{1}{2}m_{1,h}(\theta^j,\theta^j)-\frac{1}{2}m_{1,h}(\theta^{j-1},\theta^{j-1}) \right),
    \end{align*}
    and
    \begin{align*}
        m_{2,h}\left(\frac{\theta^j-\theta^{j-1}}{\tau},\theta^j\right) = \frac{1}{\tau}\left( \frac{1}{2} m_{2,h}(\theta^j-\theta^{j-1},\theta^j-\theta^{j-1})+\frac{1}{2}m_{2,h}(\theta^j,\theta^j) -\frac{1}{2}m_{2,h}(\theta^{j-1},\theta^{j-1}) \right).
    \end{align*}
    Substituting these values in \eqref{1.fullydscrterrbound1}, we have
    \begin{align*}
        \frac{1}{2\tau}\left(m_{1,h}(\theta^j,\theta^j)-m_{1,h}(\theta^{j-1},\theta^{j-1}) \right) &+ \frac{1}{2\tau}\left(m_{2,h}(\theta^j,\theta^j)-m_{2,h}(\theta^{j-1},\theta^{j-1}) \right)+ a_h(\theta^j,\theta^j) \\
        &\leq \frac{Ch^{k+1}}{\tau} \left( \|u\|_{L^2(H^{k+1}(\Omega))} + \|u_{t}\|_{L^2(H^{k+1}(\Omega))} + |u_0|_{k+1} + \tau|f^j|_{k+1} \right) \|\theta^j\|\\
        &~~+\frac{C}{\tau}\left(h^{k+1}|u^j-u^{j-1}|_{k+1} + \|\tau u_t^j-(u^j-u^{j-1})\| \right) \|\theta^j\|,
    \end{align*}
    it gives,
    \begin{align*}
        m_{1,h}(\theta^j,\theta^j)-m_{1,h}(\theta^{j-1},\theta^{j-1})&+m_{2,h}(\theta^j,\theta^j)-m_{2,h}(\theta^{j-1},\theta^{j-1})\\
        &\leq Ch^{k+1} \left( \|u\|_{L^2(H^{k+1}(\Omega))} + \|u_{t}\|_{L^2(H^{k+1}(\Omega))} + |u_0|_{k+1} + \tau|f^j|_{k+1} \right) \|\theta^j\| \\
        &~~+C\left(\|\tau u_t^j-(u^j-u^{j-1})\| +h^{k+1}|u^j-u^{j-1}|_{k+1} \right) \|\theta^j\|.
    \end{align*}
    Sum with respect to $j$ from $1$ to $n$, using $\theta^0 = \theta(\cdot,0) = 0$, it yields
    \begin{align*}
        m_{1,h}(\theta^n,\theta^n)+m_{2,h}(\theta^n,\theta^n) &\leq C \sum\limits_{j=1}^n \left(h^{k+1}\left( \|u\|_{L^2(H^{k+1}(\Omega))} + \|u_{t}\|_{L^2(H^{k+1}(\Omega))} + |u_0|_{k+1} + \tau|f^j|_{k+1} \right) \right.\\
        &~~\left. + h^{k+1}|u^j-u^{j-1}|_{k+1} + \|\tau u_t^j-(u^j-u^{j-1})\| \right) \|\theta^j\|,
    \end{align*}
    that is,
    \begin{equation}\label{1.fullydscrtbound2}
    \begin{split}
        \|\theta^n\|^2+|\theta^n|_{1}^2 &\leq C\sum\limits_{j=1}^n \left(h^{k+1}\left( \|u\|_{L^2(H^{k+1}(\Omega))} + \|u_{t}\|_{L^2(H^{k+1}(\Omega))} + |u_0|_{k+1} + \tau|f^j|_{k+1} \right) \right.\\
        &~~ \left. +h^{k+1}|u^j-u^{j-1}|_{k+1} + \|\tau u_t^j-(u^j-u^{j-1})\| \right)^2 + \frac{1}{2}\sum\limits_{j=1}^n\|\theta^j\|^2.
    \end{split}
    \end{equation}
    Using the idea of \cite{zhao2019nonconforming}, the following bounds hold:
    \begin{align}
        \tau \sum\limits_{j=1}^n h^{k+1}|f^j|_{k+1} &\leq \tau nh^{k+1} \sup\limits_{1\leq j\leq n}|f^j|_{k+1} \leq t_nh^{k+1}\|f\|_{L^{\infty}(H^{k+1}(\Omega))}, \label{1.fullydscrttaylorbnd1}\\
         \sum\limits_{j=1}^n \|\tau u_t^j-(u^j-u^{j-1})\| &\leq \tau \|u_{tt}\|_{L^1(L^2(\Omega))},
         \label{1.fullydscrttaylorbnd2}\\
         \sum\limits_{j=1}^n h^{k+1}|u^j-u^{j-1}|_{k+1} &\leq h^{k+1}\|u_{t}\|_{L^1(H^{k+1}(\Omega))}. \label{1.fullydscrttaylorbnd3}
    \end{align}
    Substituting the values from \eqref{1.fullydscrttaylorbnd1}-\eqref{1.fullydscrttaylorbnd3} into \eqref{1.fullydscrtbound2}, and it leads to 
    \begin{equation*}
        \begin{split}
            \|\theta^n\|^2+|\theta^n|_{1}^2 \leq& C \left( h^{k+1}t_n\left( \|f\|_{L^{\infty}(H^{k+1}(\Omega))} + \|u\|_{L^2(H^{k+1}(\Omega))} + \|u_{t}\|_{L^2(H^{k+1}(\Omega))} + |u_0|_{k+1} \right) \right. \\
            & \left.+ h^{k+1}\|u_{t}\|_{L^1(H^{k+1}(\Omega))} + \tau \|u_{tt}\|_{L^1(L^2(\Omega))} \right)^2 + \sum\limits_{j=1}^{n-1} \left( \|\theta^j\|^2+|\theta^j|_{1}^2 \right).
        \end{split}
    \end{equation*}
    Application of discrete Gr\"{o}nwall's inequality allows us to obtain
    \begin{equation*}\label{1.fullydscrtthetaestimate}
        \begin{split}
            \|\theta^n\|^2+|\theta^n|_{1}^2 \leq Ch^{2(k+1)} \left( \|f\|_{L^{\infty}(H^{k+1}(\Omega))}^2 + \|u\|_{L^2(H^{k+1}(\Omega))}^2 + \|u_{t}\|_{L^2(H^{k+1}(\Omega))}^2 + |u_0|_{k+1}^2\right) + C\tau^2 \|u_{tt}\|_{L^1(L^2(\Omega))}^2.
        \end{split}
    \end{equation*}
    Making use of \eqref{1.interprojerrh1} and \eqref{1.interprojerrl2} along with \eqref{1.rhobound}, we conclude that
    \begin{equation*}
        \begin{split}
            \|U^n-u^n\| + h|U^n-u^n|_{1} &\leq Ch^{k+1} \left( \|f\|_{L^{\infty}(H^{k+1}(\Omega))} + \|u\|_{L^2(H^{k+1}(\Omega))} + \|u_t\|_{L^2(H^{k+1}(\Omega))} +|u_0|_{k+1} \right) \\
            &~~+ C\tau \|u_{tt}\|_{L^1(L^2(\Omega))}.
        \end{split}
    \end{equation*}
    Hence, the proof is completed.
\end{proof}

\section{Numerical Results}
\begin{figure}[!ht]
    \centering
    \begin{minipage}{0.32\textwidth}
        \centering
        \includegraphics[width=\textwidth]{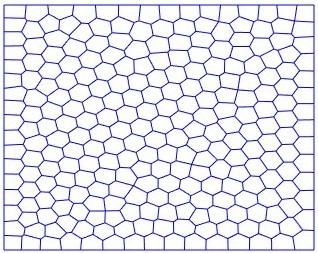}
        \caption*{\small Voronoi mesh}
    \end{minipage}
    \begin{minipage}{0.32\textwidth}
        \centering
        \includegraphics[width=\textwidth]{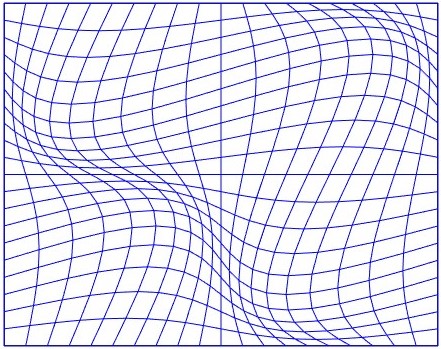}
        \caption*{\small Distorted square mesh}
    \end{minipage}
    \begin{minipage}{0.32\textwidth}
        \centering
        \includegraphics[width=\textwidth]{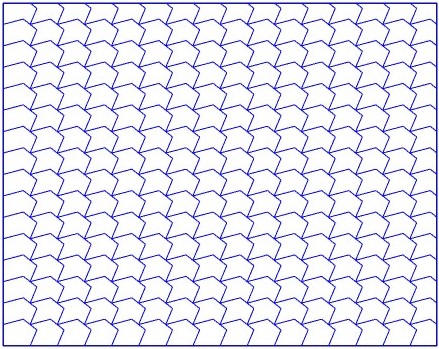}
        \caption*{\small Concave mesh}
    \end{minipage}
    \caption{Polygonal meshes}
    \label{fig:1.meshes}
\end{figure}

Several numerical experiments have been carried out to validate the theoretical results of the proposed virtual element method. In all the experiments, we have taken $\mathcal{I} = (0,1]$ and the domain is set to be $\Omega = [0,1]\times [0,1]$ discretized using Voronoi meshes $(\mathcal{V}_h)$, distorted square meshes $(\mathcal{D}_h)$ and concave polygonal meshes $(\mathcal{C}_h)$ (see Figure \ref{fig:1.meshes}). The error quantities will be computed as follows:
\begin{align*}
\begin{array}{cc}
   \mathcal{E}_{0,h}^2 = \displaystyle \sum_{K\in \mathcal{T}_h} \|u-\Pi_k^0u_h\|_K^2,  &  \mathcal{E}_{1,h}^2 = \displaystyle \sum_{K\in \mathcal{T}_h} \|\nabla u-\Pi_{k-1}^0\nabla u_h\|_K^2.
\end{array}
\end{align*}
\vspace{-1em}

\noindent \textbf{Example 1:} We consider the problem \eqref{1.mp1}-\eqref{1.mp3} in $\Omega \times \mathcal{I}$, have an analytical solution as follows:
\begin{equation*}
    u(x,y,t) = t\sin(\pi x)\sin(\pi y),
\end{equation*}
with coefficients
\begin{align*}
    \bs \mu = 
    \begin{bmatrix}
        x+y+1 & 0\\
        0 & x+y+1
    \end{bmatrix}, \ 
    \bs \eps = 
    \begin{bmatrix}
        x^2+y & 0\\
        0 & x^2+y
    \end{bmatrix}, \ 
    \bs \beta = 
    \begin{bmatrix}
        x \\
        y
    \end{bmatrix}
\end{align*}
and $\gamma = x+y$.
\begin{table}[!ht]
\begin{center}
\begin{tabular}[\textwidth]{||C{1cm}| *{5}{C{2cm}}||}
 \hline
 \textcolor{black}{$k$} & $h$ & $\mathcal{E}_{0,h}$ & EOC & $\mathcal{E}_{1,h}$ & EOC \\
 \hhline{||= = = = = =||}
 \multirow{4}{*}{1} & 3.184e-1 & 2.269e-2 & - & 5.051e-1 & - \\
 & 2.304e-1 & 1.068e-2 & 2.3294 & 3.519e-1 & 1.1171 \\
& 1.164e-1 & 2.574e-3 & 2.0829 & 1.756e-1 & 1.0179 \\
& 8.273e-2 & 1.249e-3 & 2.1199 & 1.241e-1 & 1.0169 \\
\hline
\multirow{4}{*}{2} & 3.184e-1 & 1.489e-3 & - & 6.218e-2 & - \\
& 2.304e-1 & 5.065e-4 & 3.3363 & 2.962e-2 & 2.2932 \\
& 1.164e-1 & 6.242e-5 & 3.0638 & 7.228e-3 & 2.0643 \\
& 8.273e-2 & 2.192e-5 & 3.0686 & 3.715e-3 & 1.9514 \\
\hline
\multirow{4}{*}{3} & 3.184e-1 & 1.122e-4 & - & 6.123e-3 & - \\
& 2.304e-1 & 2.383e-5 & 4.7918 & 1.864e-3 & 3.6783 \\
& 1.164e-1 & 1.408e-6 & 4.1392 & 2.203e-4 & 3.1250 \\
& 8.273e-2 & 3.349e-7 & 4.2108 & 7.579e-5 & 3.1278 \\
\hline
\end{tabular}
\caption{\label{table_eg1_voronoi} Errors and estimated orders of convergence (EOC)} for Example 1 with Voronoi mesh.
\end{center}
\end{table}

\begin{table}[!ht]
\begin{center}
\begin{tabular}[\textwidth]{||C{1cm}| *{5}{C{2cm}}||}
 \hline
 \textcolor{black}{$k$} & $h$ & $\mathcal{E}_{0,h}$ & EOC & $\mathcal{E}_{1,h}$ & EOC \\
  \hhline{||= = = = = =||}
\multirow{5}{*}{1} & 4.261e-1 & 4.377e-1 & - & 7.047e-1 & - \\
& 2.205e-1 & 1.258e-2 & 1.8925 & 3.733e-1 & 0.9643 \\
& 1.521e-1 & 5.759e-3 & 2.1052 & 2.518e-1 & 1.0613 \\
& 1.144e-1 & 3.274e-3 & 1.9826 & 1.896e-1 & 0.9956 \\
& 9.181e-2 & 2.105e-3 & 2.0055 & 1.519e-1 & 1.0054 \\
\hline
\multirow{5}{*}{2} & 4.261e-1 & 3.737e-3 & - & 1.301e-1 & - \\
& 2.205e-1 & 5.693e-4 & 2.8555 & 3.800e-2 & 1.8673 \\
& 1.521e-1 & 1.737e-4 & 3.1987 & 1.738e-2 & 2.1088 \\
& 1.144e-1 & 7.398e-5 & 2.9961 & 9.869e-3 & 1.9855 \\
& 9.181e-2 & 3.803e-5 & 3.0227 & 6.344e-3 & 2.0078 \\
\hline
\multirow{5}{*}{3} & 4.261e-1 & 3.815e-4 & - & 1.491e-2 & - \\
& 2.205e-1 & 2.948e-5 & 3.8853 & 2.186e-3 & 2.9132 \\
& 1.521e-1 & 6.033e-6 & 4.2758 & 6.607e-4 & 3.2251 \\
& 1.144e-1 & 1.928e-6 & 4.0035 & 2.802e-4 & 3.0109 \\
& 9.181e-2 & 7.930e-7 & 4.0366 & 1.437e-4 & 3.0323 \\
\hline
\end{tabular}
\caption{\label{table_eg1_distorted_sq} Errors and estimated orders of convergence for Example 1 with distorted square mesh.}
\end{center}
\end{table}
\begin{figure}[!ht]
    \centering
        \begin{minipage}{0.32\textwidth}
        \centering
        \includegraphics[width=\textwidth]{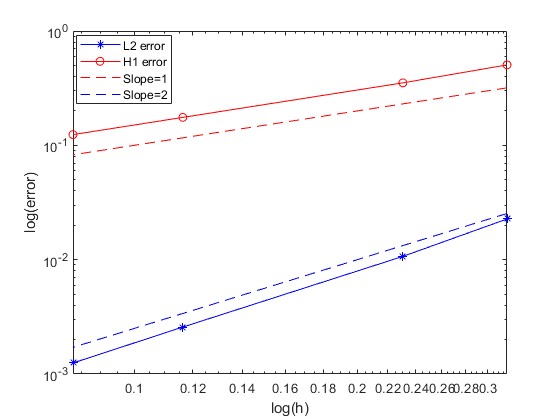}
        \caption*{\small k=1}
    \end{minipage}
    \begin{minipage}{0.32\textwidth}
        \centering
        \includegraphics[width=\textwidth]{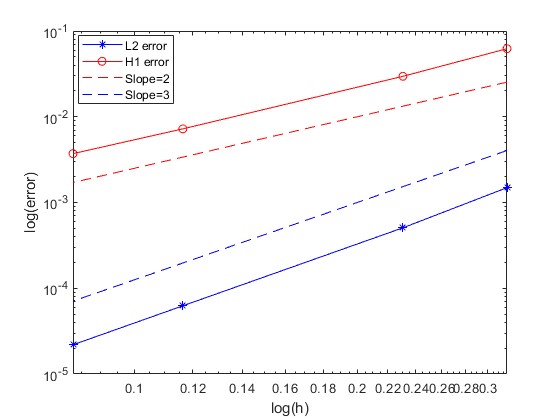}
        \caption*{\small k=2}
    \end{minipage}
    \begin{minipage}{0.32\textwidth}
        \centering
        \includegraphics[width=\textwidth]{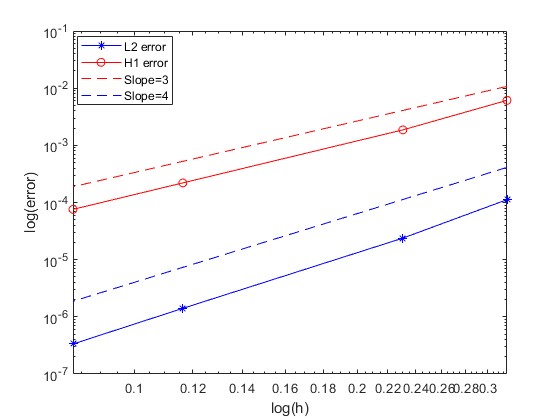}
        \caption*{\small k=3}
    \end{minipage}
    \caption{Error plots for Example $1$ on Voronoi meshes for $k = 1,\ 2$ and $3$.}
    \label{fig:error_eg1_voronoi}
    \end{figure}
   \begin{figure}[!ht]
        \begin{minipage}{0.32\textwidth}
        \centering
        \includegraphics[width=\textwidth]{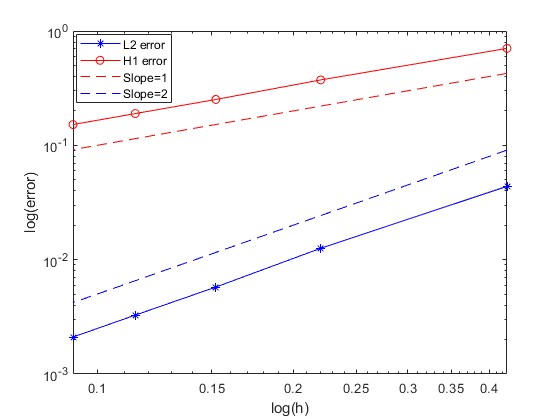}
        \caption*{\small k=1}
    \end{minipage}
    \begin{minipage}{0.32\textwidth}
        \centering
        \includegraphics[width=\textwidth]{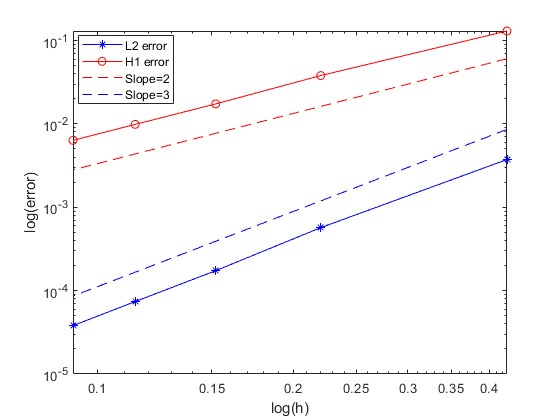}
        \caption*{\small k=2}
    \end{minipage}
    \begin{minipage}{0.32\textwidth}
        \centering
        \includegraphics[width=\textwidth]{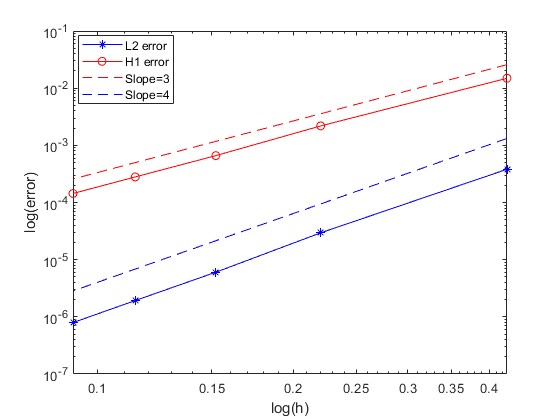}
        \caption*{\small k=3}
    \end{minipage}
    \caption{Error plots for Example $1$ on distorted square  meshes for $k = 1,\ 2$ and $3$.}
    \label{fig:error_eg1_distorted_sq}
        \end{figure}
   \begin{figure}[!ht]
        \begin{minipage}{0.32\textwidth}
        \centering
        \includegraphics[width=\textwidth]{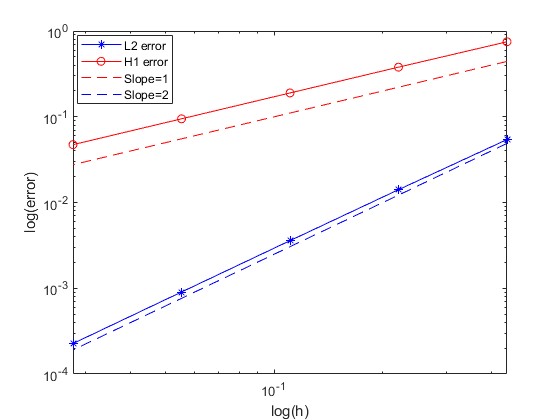}
        \caption*{\small k=1}
    \end{minipage}
    \begin{minipage}{0.32\textwidth}
        \centering
        \includegraphics[width=\textwidth]{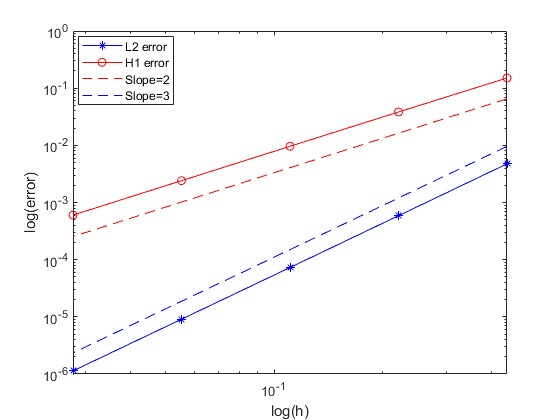}
        \caption*{\small k=2}
    \end{minipage}
    \begin{minipage}{0.32\textwidth}
        \centering
        \includegraphics[width=\textwidth]{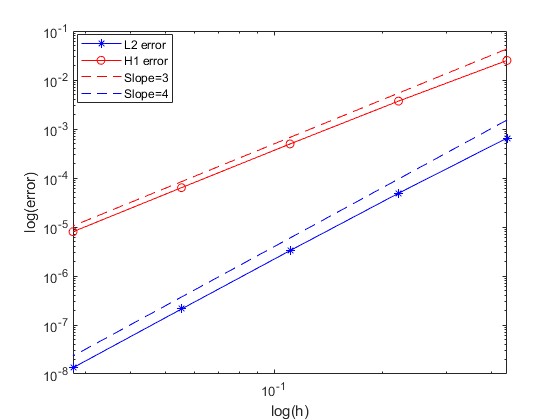}
        \caption*{\small k=3}
    \end{minipage}
    \caption{Error plots for Example $1$ on concave meshes for $k = 1,\ 2$ and $3$.}
     \label{fig:error_eg1_concave}
\end{figure}

\begin{figure}[!ht]
    \centering
        \begin{minipage}{0.32\textwidth}
        \centering
        \includegraphics[width=\textwidth]{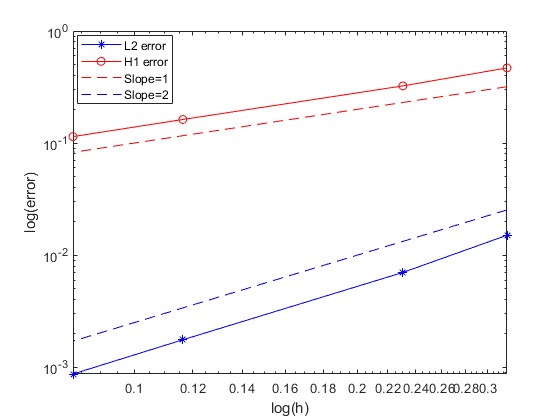}
        \caption*{\small k=1}
    \end{minipage}
    \begin{minipage}{0.32\textwidth}
        \centering
        \includegraphics[width=\textwidth]{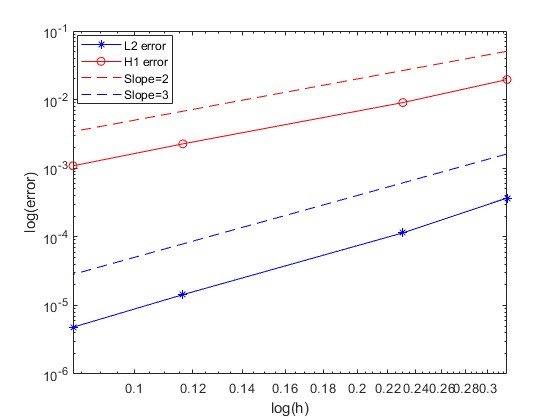}
        \caption*{\small k=2}
    \end{minipage}
    \begin{minipage}{0.32\textwidth}
        \centering
        \includegraphics[width=\textwidth]{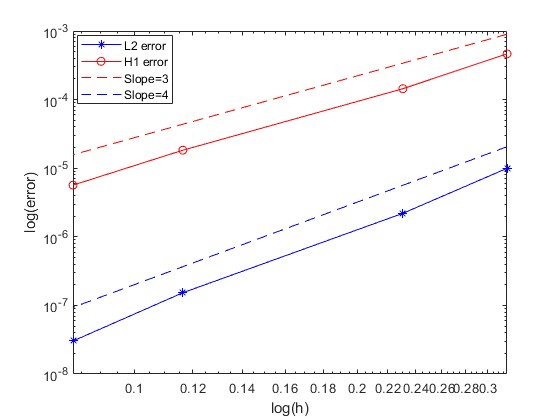}
        \caption*{\small k=3}
    \end{minipage}
    \caption{Error plots for Example $2$ on Voronoi meshes for $k = 1,\ 2$ and $3$.}
    \label{fig:error_eg2_voronoi}
\end{figure}

\begin{figure}[!ht]
        \begin{minipage}{0.32\textwidth}
        \centering
        \includegraphics[width=\textwidth]{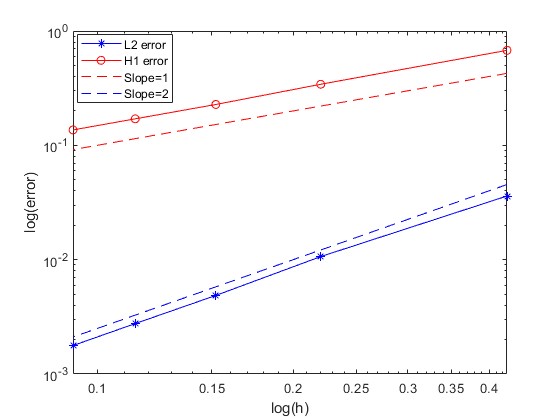}
        \caption*{\small k=1}
    \end{minipage}
    \begin{minipage}{0.32\textwidth}
        \centering
        \includegraphics[width=\textwidth]{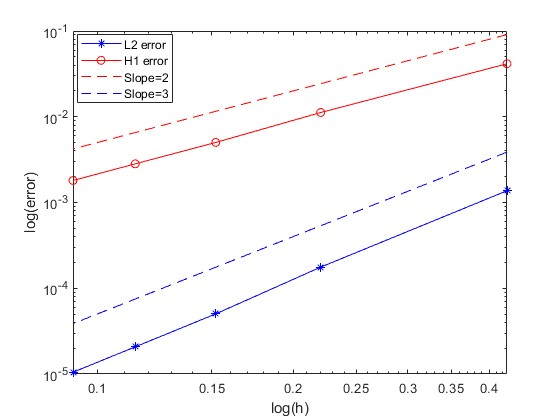}
        \caption*{\small k=2}
    \end{minipage}
    \begin{minipage}{0.32\textwidth}
        \centering
        \includegraphics[width=\textwidth]{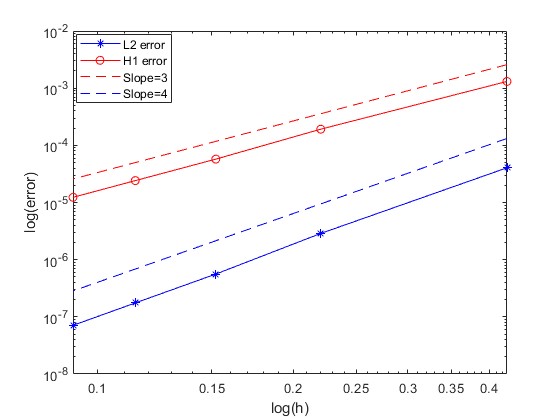}
        \caption*{\small k=3}
    \end{minipage}
    \caption{Error plots for Example $2$ on distorted square  meshes for $k = 1,\ 2$ and $3$.}
    \label{fig:error_eg2_distorted_sq}
        \end{figure}

\begin{figure}[!ht]
        \begin{minipage}{0.32\textwidth}
        \centering
        \includegraphics[width=\textwidth]{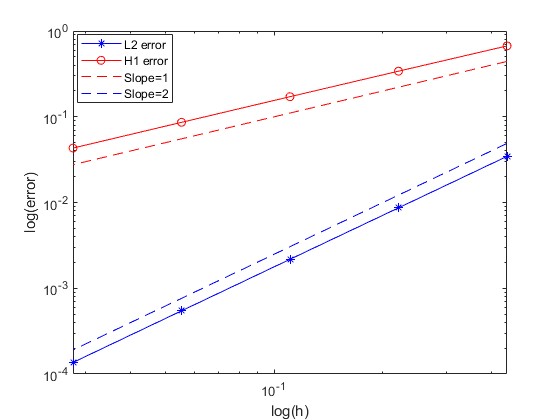}
        \caption*{\small k=1}
    \end{minipage}
    \begin{minipage}{0.32\textwidth}
        \centering
        \includegraphics[width=\textwidth]{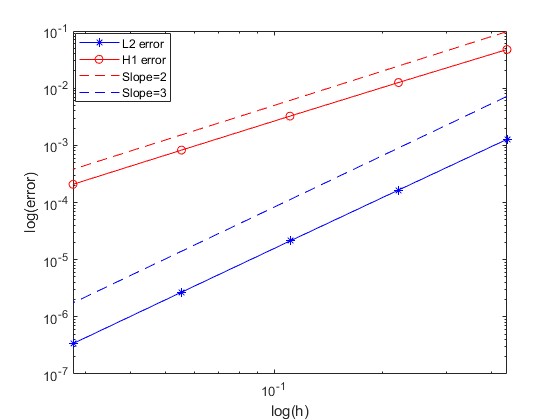}
        \caption*{\small k=2}
    \end{minipage}
    \begin{minipage}{0.32\textwidth}
        \centering
        \includegraphics[width=\textwidth]{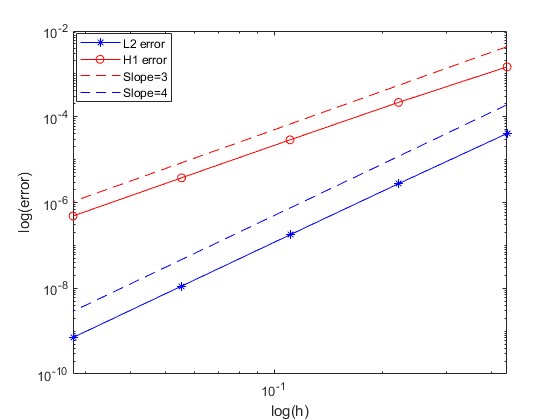}
        \caption*{\small k=3}
    \end{minipage}
    \caption{Error plots for \textcolor{black}{Example $2$} on concave meshes for $k = 1,\ 2$ and $3$.}
     \label{fig:error_eg2_concave}
\end{figure}

We can see from Figure \ref{fig:error_eg1_voronoi}, \ref{fig:error_eg1_distorted_sq} and \ref{fig:error_eg1_concave} that the numerical solution is converging in both $L^2-$norm and energy norm with the rate of $\textit{O}(h^{k+1})$ and $\textit{O}(h^{k})$, respectively, for Voronoi, distorted square and concave mesh, respectively. Here, in each case we have chosen $\tau = 10^{-3}$, mesh size $h$ and the corresponding error values \textcolor{black}{along with the estimated order of convergence (EOC)} are given in Table \ref{table_eg1_voronoi}, \ref{table_eg1_distorted_sq} and \ref{table_eg1_concave_sq}.
\vskip5pt
\noindent \textbf{Example 2:} We now solve a convection-dominated-diffusion problem in $\Omega \times \mathcal{I}$ with non-homogeneous boundary condition, the coefficients $\bs \mu = I, \bs \eps = 10^{-6}I, \beta =   
    \begin{bmatrix}
        10 \\
        10
    \end{bmatrix}$ and $\gamma = 1$. We set the source term and the boundary condition according to the exact solution:
    \begin{equation*}
        u(x,y,t) = te^{x+y}.
    \end{equation*}
    The errors and estimated order of convergence in both $L^2-$norm and energy norm have been shown in Table \ref{table_eg2_voronoi} for Voronoi mesh, Table \ref{table_eg2_distorted_sq} for distorted square mesh and Table \ref{table_eg2_concave_sq} for concave mesh. Figure \ref{fig:error_eg2_voronoi}, \ref{fig:error_eg2_distorted_sq} and \ref{fig:error_eg2_concave} presents the log-log plot of the errors with Voronoi, distorted square, and concave mesh, respectively.

\begin{figure}[!ht]
    \centering
    \begin{minipage}{0.45\textwidth}
        \centering
        \includegraphics[width=\textwidth]{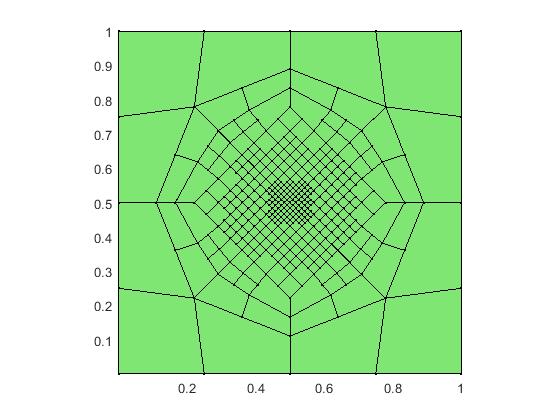}
    \end{minipage}
    \begin{minipage}{0.45\textwidth}
        \centering
        \includegraphics[width=\textwidth]{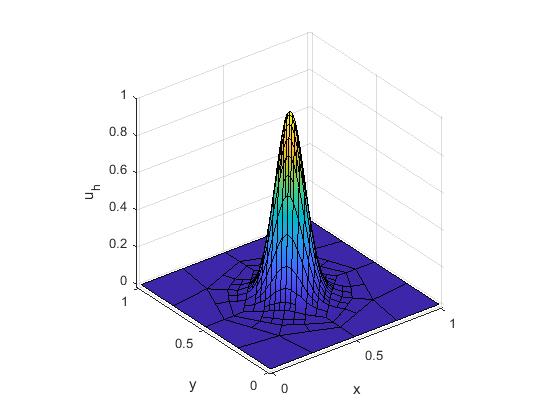}
    \end{minipage}
    \caption{The left image shows the domain discretized by adaptively refined mesh and right image is the approximate solution on it.}
    \label{fig:1.adaptive_mesh}
\end{figure}

\begin{figure}[!ht]
    \centering
    \begin{minipage}{0.45\textwidth}
        \centering
        \includegraphics[width=\textwidth]{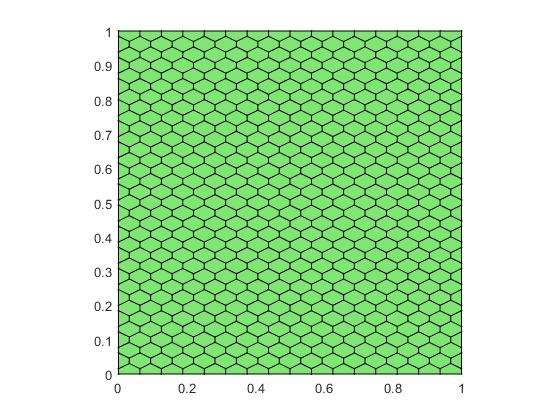}
    \end{minipage}
    \begin{minipage}{0.45\textwidth}
        \centering
        \includegraphics[width=\textwidth]{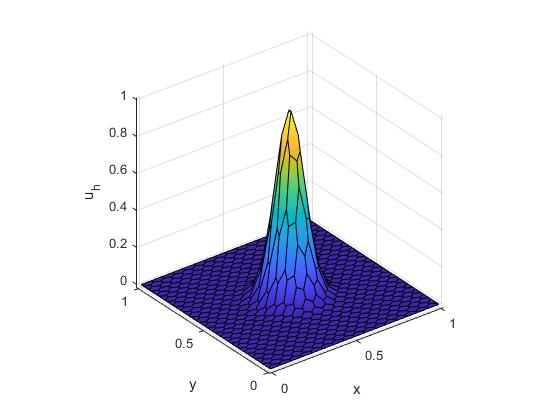}
    \end{minipage}
    \caption{The left image shows the domain discretized by uniformly refined mesh, and the right image is the approximate solution on it.}
    \label{fig:1.uniform_mesh}
\end{figure}

\begin{figure}[!ht]
    \centering
    \includegraphics[width=0.5\linewidth]{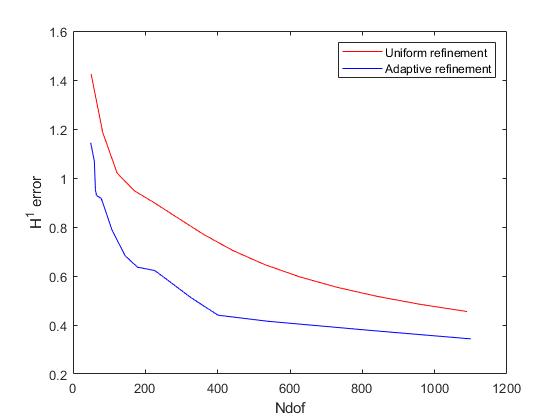}
    \caption{Comparison of errors with uniform and adaptive mesh refinement.}
    \label{fig:1.comparison}
\end{figure}

\vskip5pt
\noindent \textbf{Example 3:}  Let us consider the problem \eqref{1.mp1}-\eqref{1.mp3} in $\Omega \times \mathcal{I}$ with coefficients $\bs \mu = I, \bs \eps = I, \beta = 
    \begin{bmatrix}
        1 \\
        1
    \end{bmatrix}$ and $\gamma = 1$ have an analytical solution as:
    \begin{equation*}
        u(x,y,t) = te^{-100((x-1/2)^2+(y-1/2)^2)},
    \end{equation*}
    having a sharp peek at the point $(\frac{1}{2},\frac{1}{2})$. Therefore, adaptive mesh refinement is a better way to go than doing it uniformly. The Figure \ref{fig:1.adaptive_mesh} 
    is displaying the adaptively refined mesh and the approximate solution on it while the findings for uniformly refined mesh are shown in Figure \ref{fig:1.uniform_mesh}. Figure \ref{fig:1.adaptive_mesh} illustrates that the adaptively refined mesh has hanging nodes, which can be handled with VEM, whereas further elements would have needed to refine in FEM. In Figure \ref{fig:1.comparison}, we are showing that for the same number of degrees of freedom, the solution obtained by adaptively refined mesh outperforms the solution obtained by uniformly refined mesh in terms of accuracy.

\begin{table}[!ht]
\begin{center}
\begin{tabular}[\textwidth]{||C{1cm}| *{5}{C{2cm}}||}
 \hline
 \textcolor{black}{$k$} & $h$ & $\mathcal{E}_{0,h}$ & EOC & $\mathcal{E}_{1,h}$ & EOC \\
\hhline{||= = = = = =||}
\multirow{5}{*}{1} & 4.419e-1 & 50420e-2 & - & 7.487e-1 & - \\
& 2.209e-1 & 1.411e-2 & 1.9412 & 3.788e-1 & 0.9830 \\
& 1.105e-1 & 3.589e-3 & 1.9754 & 1.894e-1 & 0.9996 \\
& 5.524e-2 & 9.039e-4 & 1.9895 & 9.462e-2 & 1.0017 \\
& 2.762e-2 & 2.2167e-4 & 1.9954 & 4.727e-2 & 1.0012 \\
\hline
\multirow{5}{*}{2} & 4.419e-1 & 4.747e-3 & - & 1.507e-1 & - \\
& 2.209e-1 & 5.875e-4 & 3.0142 & 3.831e-2 & 1.9757 \\
& 1.105e-1 & 7.296e-5 & 3.0095 & 9.621e-3 & 1.9937 \\
& 5.524e-2 & 9.087e-6 & 3.0051 & 2.408e-3 & 1.9981 \\
& 2.762e-2 & 1.134e-6 & 3.0026 & 6.023e-4 & 1.9994 \\
\hline
\multirow{5}{*}{3} & 4.419e-1 & 6.523e-4 & - & 2.508e-2 & - \\
& 2.209e-1 & 4.881e-5 & 3.7403 & 3.706e-3 & 2.7585 \\
& 1.105e-1 & 3.295e-6 & 3.8889 & 4.953e-4 & 2.9037 \\
& 5.524e-2 & 2.121e-7 & 3.9571 & 6.347e-5 & 2.9641 \\
& 2.762e-2 & 1.342e-8 & 3.9823 & 8.013e-6 & 2.9856 \\
\hline
\end{tabular}
\caption{\label{table_eg1_concave_sq} Errors and estimated orders of convergence for Example 1 with concave mesh.}
\end{center}
\end{table}

\begin{table}[!ht]
\begin{center}
\begin{tabular}[\textwidth]{||C{1cm}| *{5}{C{2cm}}||}
 \hline
 \textcolor{black}{$k$} & $h$ & $\mathcal{E}_{0,h}$ & EOC & $\mathcal{E}_{1,h}$ & EOC \\
 \hhline{||= = = = = =||}
\multirow{5}{*}{1} & 4.261e-1 & 3.603e-2 & - & 6.774e-1 & - \\
& 2.205e-1 & 1.064e-2 & 1.8509 & 3.421e-1 & 1.0366 \\
& 1.521e-1 & 4.864e-3 & 2.1095 & 2.277e-1 & 1.0976 \\
& 1.144e-1 & 2.761e-3 & 1.9874 & 1.706e-1 & 1.0136 \\
& 9.181e-2 & 1.774e-3 & 2.0090 & 1.364e-1 & 1.0165 \\
\hline
\multirow{5}{*}{2} & 4.261e-1 & 1.371e-3 & - & 4.147e-2 & - \\
& 2.205e-1 & 1.764e-4 & 3.1123 & 1.121e-2 & 1.9854 \\
& 1.521e-1 & 5.027e-5 & 3.3826 & 5.017e-3 & 2.1666 \\
& 1.144e-1 & 2.074e-5 & 3.1073 & 2.823e-3 & 2.0187 \\
& 9.181e-2 & 1.049e-5 & 3.0986 & 1.805e-3 & 2.0302 \\
\hline
\multirow{5}{*}{3} & 4.261e-1 & 4.048e-5 & - & 1.312e-3 & - \\
& 2.205e-1 & 2.874e-6 & 4.0140 & 1.917e-4 & 2.9185 \\
& 1.521e-1 & 5.600e-7 & 4.4074 & 5.736e-5 & 3.2512 \\
& 1.144e-1 & 1.748e-7 & 4.0868 & 2.414e-5 & 3.0377 \\
& 9.181e-2 & 7.097e-8 & 4.0955 & 1.232e-5 & 3.0562 \\
\hline
\end{tabular}
\caption{\label{table_eg2_distorted_sq} Errors and estimated orders of convergence for Example 2 with distorted square mesh.}
\end{center}
\end{table}

\begin{table}[!ht]
\begin{center}
\begin{tabular}[\textwidth]{||C{1cm}| *{5}{C{2cm}}||}
 \hline
 \textcolor{black}{$k$} & $h$ & $\mathcal{E}_{0,h}$ & EOC & $\mathcal{E}_{1,h}$ & EOC \\
 \hhline{||= = = = = =||}
\multirow{4}{*}{1} & 3.184e-1 & 1.506e-2 & - & 4.683e-1 & - \\
& 2.304e-1 & 7.012e-3 & 2.3635 & 3.237e-1 & 1.1423 \\
& 1.164e-1 & 1.759e-3 & 2.0241 & 1.624e-1 & 1.0097 \\
& 8.273e-2 & 8.704e-4 & 2.0618 & 1.1423e-1 & 1.0309 \\
\hline
\multirow{4}{*}{2} & 3.184e-1 & 3.677e-4 & - & 1.959e-2 & - \\
& 2.304e-1 & 1.142e-4 & 3.6174 & 9.047e-3 & 2.3905 \\
& 1.164e-1 & 1.434e-5 & 3.0361 & 2.263e-3 & 2.0278 \\
& 8.273e-2 & 4.807e-6 & 3.2036 & 1.084e-3 & 2.1567 \\
\hline
\multirow{4}{*}{3} & 3.184e-1 & 9.959e-6 & - & 4.639e-4 & - \\
& 2.304e-1 & 2.200e-6 & 4.6695 & 1.441e-4 & 3.6152 \\
& 1.164e-1 & 1.529e-7 & 3.9015 & 1.829e-5 & 3.0205 \\
& 8.273e-2 & 3.029e-8 & 4.7476 & 5.663e-6 & 3.4381 \\
\hline
\end{tabular}
\caption{\label{table_eg2_voronoi} Errors and estimated orders of convergence for Example 2 with Voronoi mesh.}
\end{center}
\end{table}
 
\begin{table}[!ht]
\begin{center}
\begin{tabular}[\textwidth]{||C{1cm}| *{5}{C{2cm}}||}
 \hline
 \textcolor{black}{$k$} & $h$ & $\mathcal{E}_{0,h}$ & EOC & $\mathcal{E}_{1,h}$ & EOC \\
 \hhline{||= = = = = =||}
\multirow{5}{*}{1} & 4.419e-1 & 3.451e-2 & - & 6.700e-1 & - \\
& 2.209e-1 & 8.662e-3 & 1.9941 & 3.401e-1 & 0.9784 \\
& 1.105e-1 & 2.176e-3 & 1.9928 & 1.714e-1 & 0.9938 \\
& 5.524e-2 & 5.452e-4 & 1.9969 & 8.607e-2 & 0.9938 \\
& 2.762e-2 & 1.364e-4 & 1.9988 & 4.313e-2 & 0.9968 \\
\hline
\multirow{5}{*}{2} & 4.419e-1 & 1.275e-3 & - & 4.742e-2 & - \\
& 2.209e-1 & 1.651e-4 & 2.9482 & 1.256e-2 & 1.9170 \\
& 1.105e-1 & 2.123e-5 & 2.9597 & 3.244e-3 & 1.9524 \\
& 5.524e-2 & 2.698e-6 & 2.9761 & 8.253e-4 & 1.9749 \\
& 2.762e-2 & 3.402e-7 & 2.9871 & 2.082e-4 & 1.9871 \\
\hline
\multirow{5}{*}{3} & 4.419e-1 & 4.061e-5 & - & 1.456e-3 & - \\
& 2.209e-1 & 2.747e-6 & 3.8862 & 2.161e-4 & 2.7528 \\
& 1.105e-1 & 1.764e-7 & 3.9604 & 2.894e-5 & 2.9006 \\
& 5.524e-2 & 1.120e-8 & 3.9773 & 3.761e-6 & 2.9442 \\
& 2.762e-2 & 7.091e-10 & 3.9816 & 4.841e-7 & 2.9573 \\
\hline
\end{tabular}
\caption{\label{table_eg2_concave_sq} Errors and estimated orders of convergence for Example 2 with concave mesh.}
\end{center}
\end{table}

\section{Conclusion}
In this paper, we have presented a VEM for the Sobolev equation with the convection term and with variable coefficients. Our newly introduced intermediate projection allowed us to achieve the best possible order of convergence in both $L^2$-norm and the energy norm. Additionally, the optimal convergence was determined by employing Euler's backward method to discretize the temporal variable. Then, we conducted numerous numerical experiments that demonstrated complete agreement with the theoretical results. Our future scope is to study the Sobolev equation using non-conforming VEM.\\

\noindent \textcolor{black}{\textbf{Acknowledgements} The authors thank the anonymous reviewer for the valuable comments and suggestions which led to numerous improvements to the manuscript. The second author would like to thank Anusandhan National Research Foundation (ANRF), a statutory body of the Department of Science and Technology (DST), for supporting this work through the core research grant CRG/2021/002410.}\\

\noindent \textbf{Conflict of interest} The authors have no conflict of interest in this work.

\bibliographystyle{siamplain}
\bibliography{ref}

\begin{thebibliography}{10}

\bibitem{ahmad2013equivalent}
{\sc B.~Ahmad, A.~Alsaedi, F.~Brezzi, L.~D. Marini, and A.~Russo}, {\em
  Equivalent projectors for virtual element methods}, Computers \& Mathematics
  with Applications, 66 (2013), pp.~376--391.

\bibitem{anaya2020virtual}
{\sc V.~Anaya, M.~Bendahmane, D.~Mora, and M.~Sep{\'u}lveda}, {\em A virtual
  element method for a nonlocal {F}itzhugh--{N}agumo model of cardiac
  electrophysiology}, IMA Journal of Numerical Analysis, 40 (2020),
  pp.~1544--1576.

\bibitem{antonietti2014stream}
{\sc P.~F. Antonietti, L.~B. Da~Veiga, D.~Mora, and M.~Verani}, {\em A stream
  virtual element formulation of the {S}tokes problem on polygonal meshes},
  SIAM Journal on Numerical Analysis, 52 (2014), pp.~386--404.

\bibitem{antonietti2023virtual}
{\sc P.~F. Antonietti, G.~Vacca, and M.~Verani}, {\em Virtual element method
  for the {N}avier--{S}tokes equation coupled with the heat equation}, IMA
  Journal of Numerical Analysis, 43 (2023), pp.~3396--3429.

\bibitem{arrutselvi2021virtual}
{\sc M.~Arrutselvi and E.~Natarajan}, {\em Virtual element method for nonlinear
  time-dependent convection-diffusion-reaction equation}, Computational
  Mathematics and Modeling, 32 (2021), pp.~376--386.

\bibitem{barenblatt1960basic}
{\sc G.~I. Barenblatt, I.~P. Zheltov, and I.~Kochina}, {\em Basic concepts in
  the theory of seepage of homogeneous liquids in fissured rocks [strata]},
  Journal of applied mathematics and mechanics, 24 (1960), pp.~1286--1303.

\bibitem{beirao2013basic}
{\sc L.~Beir{\~a}o~da Veiga, F.~Brezzi, A.~Cangiani, G.~Manzini, L.~D. Marini,
  and A.~Russo}, {\em Basic principles of virtual element methods},
  Mathematical Models and Methods in Applied Sciences, 23 (2013), pp.~199--214.

\bibitem{beirao2016virtual}
{\sc L.~Beir{\~a}o~da Veiga, F.~Brezzi, L.~D. Marini, and A.~Russo}, {\em
  Virtual element method for general second-order elliptic problems on
  polygonal meshes}, Mathematical Models and Methods in Applied Sciences, 26
  (2016), pp.~729--750.

\bibitem{benedetto2016order}
{\sc M.~F. Benedetto, S.~Berrone, A.~Borio, S.~Pieraccini, and S.~Scialo}, {\em
  Order preserving {SUPG} stabilization for the virtual element formulation of
  advection--diffusion problems}, Computer Methods in Applied Mechanics and
  Engineering, 311 (2016), pp.~18--40.

\bibitem{brenner2017some}
{\sc S.~C. Brenner, Q.~Guan, and L.-Y. Sung}, {\em Some estimates for virtual
  element methods}, Computational Methods in Applied Mathematics, 17 (2017),
  pp.~553--574.

\bibitem{cangiani2017conforming}
{\sc A.~Cangiani, G.~Manzini, and O.~J. Sutton}, {\em Conforming and
  nonconforming virtual element methods for elliptic problems}, IMA Journal of
  Numerical Analysis, 37 (2017), pp.~1317--1354.

\bibitem{chen2019two}
{\sc C.~Chen, K.~Li, Y.~Chen, and Y.~Huang}, {\em Two-grid finite element
  methods combined with {C}rank-{N}icolson scheme for nonlinear {S}obolev
  equations}, Advances in Computational Mathematics, 45 (2019), pp.~611--630.

\bibitem{da2013virtual}
{\sc L.~B. Da~Veiga, F.~Brezzi, and L.~D. Marini}, {\em Virtual elements for
  linear elasticity problems}, SIAM Journal on Numerical Analysis, 51 (2013),
  pp.~794--812.

\bibitem{da2021supg}
{\sc L.~B. da~Veiga, F.~Dassi, C.~Lovadina, and G.~Vacca}, {\em
  {SUPG}-stabilized virtual elements for diffusion-convection problems: a
  robustness analysis}, ESAIM: Mathematical Modelling and Numerical Analysis,
  55 (2021), pp.~2233--2258.

\bibitem{dongyang2016unconditional}
{\sc S.~Dongyang, Y.~Fengna, and W.~Junjun}, {\em Unconditional
  superconvergence analysis of a new mixed finite element method for nonlinear
  {S}obolev equation}, Applied Mathematics and Computation, 274 (2016),
  pp.~182--194.

\bibitem{ewing1978time}
{\sc R.~E. Ewing}, {\em Time-stepping {G}alerkin methods for nonlinear
  {S}obolev partial differential equations}, SIAM Journal on Numerical
  Analysis, 15 (1978), pp.~1125--1150.

\bibitem{gao2017weak}
{\sc F.~Gao, J.~Cui, and G.~Zhao}, {\em Weak {G}alerkin finite element methods
  for {S}obolev equation}, Journal of Computational and Applied Mathematics,
  317 (2017), pp.~188--202.

\bibitem{gao2009local}
{\sc F.~Gao, J.~Qiu, and Q.~Zhang}, {\em Local discontinuous {G}alerkin finite
  element method and error estimates for one class of {S}obolev equation},
  Journal of Scientific Computing, 41 (2009), pp.~436--460.

\bibitem{gao2009split}
{\sc F.~Gao and H.~Rui}, {\em A split least-squares characteristic mixed finite
  element method for {S}obolev equations with convection term}, Mathematics and
  Computers in Simulation, 80 (2009), pp.~341--351.

\bibitem{kumar2023stabilizer}
{\sc N.~Kumar and B.~Deka}, {\em A stabilizer free weak {G}alerkin finite
  element method for second-order {S}obolev equation}, Numerical Methods for
  Partial Differential Equations, 39 (2023), pp.~2115--2140.

\bibitem{kumar2024numerical}
{\sc S.~Kumar, D.~Mora, R.~Ruiz-Baier, and N.~Verma}, {\em Numerical solution
  of the {B}iot/elasticity interface problem using virtual element methods},
  Journal of Scientific Computing, 98 (2024), p.~53.

\bibitem{larson2013finite}
{\sc M.~G. Larson and F.~Bengzon}, {\em The finite element method: theory,
  implementation, and applications}, vol.~10, Springer Science \& Business
  Media, 2013.

\bibitem{li2019expanded}
{\sc N.~Li, P.~Lin, and F.~Gao}, {\em An expanded mixed finite element method
  for two-dimensional {S}obolev equations}, Journal of Computational and
  Applied Mathematics, 348 (2019), pp.~342--355.

\bibitem{mendina2012sensitivity}
{\sc M.~Mendina and R.~Terra}, {\em Sensitivity of simulated convection to soil
  moisture in a region in central amazon}, Atm{\'o}sfera, 25 (2012),
  pp.~269--293.

\bibitem{pradhan2024optimal}
{\sc G.~Pradhan and B.~Deka}, {\em Optimal convergence analysis of the virtual
  element methods for second-order {S}obolev equations with variable
  coefficients on polygonal meshes}, Journal of Applied Mathematics and
  Computing, 70 (2024), pp.~2313--2341.

\bibitem{showalter1975sobolev}
{\sc R.~Showalter}, {\em The {S}obolev equation, ii}, Applicable Analysis, 5
  (1975), pp.~81--99.

\bibitem{ting1974cooling}
{\sc T.~W. Ting}, {\em A cooling process according to two-temperature theory of
  heat conduction}, Journal of Mathematical Analysis and Applications, 45
  (1974), pp.~23--31.

\bibitem{vacca2015virtual}
{\sc G.~Vacca and L.~Beir{\~a}o~da Veiga}, {\em Virtual element methods for
  parabolic problems on polygonal meshes}, Numerical Methods for Partial
  Differential Equations, 31 (2015), pp.~2110--2134.

\bibitem{wang2024influence}
{\sc G.~Wang, R.~Fu, Y.~Zhuang, P.~A. Dirmeyer, J.~A. Santanello, G.~Wang,
  K.~Yang, and K.~McColl}, {\em Influence of lower-tropospheric moisture on
  local soil moisture--precipitation feedback over the {US} {S}outhern {G}reat
  {P}lains}, Atmospheric Chemistry and Physics, 24 (2024), pp.~3857--3868.

\bibitem{xie2022hybrid}
{\sc C.-M. Xie, M.-F. Feng, and Y.~Luo}, {\em A hybrid high-order method for
  the {S}obolev equation}, Applied Numerical Mathematics, 178 (2022),
  pp.~84--97.

\bibitem{xu2022conforming}
{\sc Y.~Xu, Z.~Zhou, and J.~Zhao}, {\em Conforming virtual element methods for
  {S}obolev equations}, Journal of Scientific Computing, 93 (2022), p.~32.

\bibitem{yadav2024conforming}
{\sc S.~Yadav, M.~Suthar, and S.~Kumar}, {\em A conforming virtual element
  method for parabolic integro-differential equations}, Computational Methods
  in Applied Mathematics, 24 (2024), pp.~1001--1019.

\bibitem{zhang2023virtual}
{\sc B.~Zhang, J.~Zhao, and S.~Chen}, {\em Virtual element method for the
  {S}obolev equations}, Mathematical Methods in the Applied Sciences, 46
  (2023), pp.~1266--1281.

\bibitem{zhao2019nonconforming}
{\sc J.~Zhao, B.~Zhang, and X.~Zhu}, {\em The nonconforming virtual element
  method for parabolic problems}, Applied Numerical Mathematics, 143 (2019),
  pp.~97--111.

\end{thebibliography}
\end{document}